\RequirePackage{fix-cm}
\documentclass{svjour3}                     
\smartqed  
\usepackage{graphicx,color}
\usepackage{amsmath,amssymb}
\usepackage{float,booktabs}
\usepackage{psfrag,multirow,tikz}
\usepackage{algorithm}
\usepackage{algorithmic}
\begin{document}

\title{Adaptive two- and three-dimensional multiresolution computations of resistive magnetohydrodynamics
}

\titlerunning{Adaptive multiresolution computations of resistive magnetohydrodynamics}        

\author{Anna Karina Fontes Gomes \and \\
		Margarete Oliveira Domingues \and \\
		Odim Mendes \and \\
		Kai Schneider
}


\institute{Anna Karina Fontes Gomes \at
              Federal Institute of S\~ao Paulo, Cubat\~ao, S\~ao Paulo, Brazil \\
              \email{anna.gomes@ifsp.edu.br}           
           \and
           Margarete Oliveira Domingues \at
              National Institute for Space Research, S\~ao Jos\'e dos Campos, S\~ao Paulo, Brasil
              \email{margarete.domingues@inpe.br}
           \and
           Odim Mendes \at
              National Institute for Space Research, S\~ao Jos\'e dos Campos, S\~ao Paulo, Brasil
              \email{odim.mendes@inpe.br}
             \and
           Kai Schneider \at
              Aix-Marseille Universit\'e, CNRS, Centrale Marseille, I2M, Marseille, France
              \email{kai.schneider@univ-amu.fr}
}

\date{Preprint accepted at Advances in Computational Mathematics}

\maketitle

\begin{abstract}
Fully adaptive computations of the resistive magnetohydrodynamic (MHD) equations are presented in two and three space dimensions using a finite volume discretization on locally refined dyadic grids.
Divergence cleaning is used to control the incompressibility constraint of the magnetic field. 
For automatic grid adaptation a cell-averaged multiresolution analysis is applied which guarantees the precision of the adaptive computations, while reducing CPU time and memory requirements.
Implementation issues of the open source code CARMEN-MHD are discussed. To illustrate its precision and efficiency different benchmark computations including shock-cloud interaction and magnetic reconnection are presented.
\end{abstract}

\section{Introduction}
\label{intro}
The constant need for understanding the nonlinear dynamics of different phenomena encountered in our daily life, which are typically governed by nonlinear partial differential equations (PDEs), calls for robust and efficient numerical methods to perform high fidelity numerical simulations.
Many complex processes, necessitate high resolution computations to represent efficiently the dynamics of a given multiscale problem~\cite{Multiscale:2014}. 
Increasing the resolution of the computational mesh directly impacts on the computational cost, which thus increases at best linearly, and thus can lead in many cases to computationally prohibitive simulations. 
In this context, dynamically adaptive multiscale methods play a prominent and important role, since their purpose is to adapt the computational mesh to the local structures present in the numerical solution, while preserving the accuracy of the adaptive computations. In particular, many phenomena in space physics can benefit from such adaptive approaches due to their intrinsic multiscale characteristics.
The presence of multiple time and space scales appears to be optimal for adaptive methodologies and highly compressed data representation, see e.g. \cite{kolomenskiy2018data}.

Different adaptive discretization schemes for magnetohydrodynamic (MHD) simulations have been proposed in the literature, an exhaustive review is beyond the scope of the present work. 
In the following we briefly describe related and competitive adaptive MHD approaches.
In~\cite{fambri2017spacetime} a space-time adaptive method using high order discontinuous Galerkin discretizations (ADER-DG) is proposed and among others applied  to viscous and resistive MHD in two and three space dimensions. This work is based on the ADER-DG schemes developed for hyperbolic conservation laws~\cite{zanotti2015spacetime}.
A parallel MHD code, NIRVANA, using adaptive mesh refinement with block-structure and domain decomposition is presented in~\cite{ziegler2008nirvana}.
A combination of adaptive mesh refinement and central weighted essentially non-oscillatory schemes has been put forward in~\cite{kleimann2004three}. The resulting third order accurate scheme has been applied to highly super-Alfvenic plasmas and to stiff Sedov-type explosion problems.
A robust second order, shock-capturing numerical scheme for multidimensional special relativistic magnetohydrodynamics can be found in~\cite{vanderholst2008multi}, again in the framework of adaptive mesh refinement and using a finite volume setting. Applications to relativistic MHD Riemann problems for which exact solutions are known, are shown to be successfully recovered.

{A detailed discussion on space weather forecasting can be found in~\cite{Toth:2012}, focusing on a publicly available Space Weather Modeling Framework (SWMF).
The foundations are a Block-Adaptive Tree Solarwind Roe-type Upwind Scheme (BATS-R-US) code that can solve various forms of the  MHD equations, including Hall, semi-relativistic, multi-species and multi-fluid MHD, anisotropic pressure, radiative transport and heat conduction. A block-adaptive mesh in Cartesian and generalized coordinates is used together with load balancing and message passing for one, two and three-dimensional problems. 
Time-stepping of SWMF can be either explicit, semi-implicit or fully implicit, depending on the application featuring likewise local time-stepping.
The current status of MHD simulations for space weather is reviewed in the recent book of Feng \cite{feng:2020}, including AMR and data driven MHD modeling within the framework of cell-centered finite volume methods.}

{The review of Jardin~\cite{Jardin:2012} discusses the importance of implicit algorithms in the context of magnetically confined fusion plasma using the MHD description. A combination of implicit solvers with highly accurate spatial discretizations and anisotropic thermal conduction is shown to allow predicting accurately fusion experiments for realistic physical parameters in realistic toroidal geometries.}

{A Lagrangian parallel MHD code, GRADSPMHD,  based on the Smooth Particle Hydrodynamics (SPH) formalism is introduced in~\cite{Vanaverbeke:2014}. A mixed hyperbolic-parabolic correction scheme is used for satisfying the divergence constraint on the magnetic field and a tree-based code for finding the neighbors.
For validation classical benchmarks were computed, including the magneto-rotational instability and simulations of magnetized accretion disks. The performance of the code on a parallel supercomputer with distributed memory architecture is likewise assessed.}

{Meshless finite-volume Lagrangian methods for hydrodynamics have been extended for ideal MHD in \cite{hopkins:2016} using a divergence cleaning scheme. Benchmark computations 
illustrate that the developed code GIZMO is competitive with adaptive mesh refinement (AMR) techniques.
Compared to SPH these methods allow sharp shock-capturing, reduced noise, divergence errors, and diffusion. However the convergence of the method is found to be problem dependent.
}

Here we present an alternative to AMR, which is meanwhile a standard approach for solving PDEs on adaptive grids \cite{berger1984adaptive,bergercollela1989}.
We propose using multiresolution (MR) analysis for introducing adaptivity in MHD simulations. 
MR is based on the idea that a data set (i.e. the solution of the PDE) can be represented at different refinement levels, according to its local regularity. 
A detailed comparison of MR and AMR approaches has been carried out in \cite{Deiterdingetal:2009,deiterding2016comparison} for compressible Euler equations. There we found that the MR method yields a better memory compression than AMR together with improved convergence.

In particular, we focus on the adaptive multiresolution for cell averages, firstly introduced by 
Harten \cite{HARTI1993153,harten1994adaptive,Harten:1995,Harten:1996} in one dimension, which is directly related to biorthogonal wavelets. The wavelet coefficients provide the information about the regularity of the data, which are used to adapt the computational mesh to the problem of interest. After Harten's seminal work, numerous publications contributed to the development of this approach in a way that the local regularity of the solution is detected \cite{Cohen2003,Kaibara2001,Mueller:2003,Roussel:2003}. Later, the adaptive MR for cell averages was extended for two and three dimensions \cite{bihari1997multiresolution,Roussel:2003}, making it possible to apply this methodology to different problems of practical interest. 
%
{In the context of point-value MR an  adaptive solver for the two-dimensional compressible Euler equations was proposed in~\cite{Chiaavassa2001}. For a detailed review on adaptive MR and {wavelet} methods for conservation laws and applications in computational fluid dynamics {solving the Navier--Stokes equations}, we refer to \cite{Mueller:2003,schneider2010wavelet,DGRSESAIM:2011}}.

In this work, we combine the finite volume method with an adaptive MR approach to solve numerically the resistive magnetohydrodynamic equations, as discussed in \cite{domingues2013extended,Gomesetal:2015} for the ideal MHD. Magnetohydrodynamics describes the behavior of a macroscopic electrically conducting fluid, which can be used to model the dynamics of space plasma~\cite{goedbloed2004principles}. The MHD model is characterized by a set of nonlinear evolutionary partial differential equations, presenting in some cases strong discontinuities of the solution. In order to solve these equations, we use robust numerical schemes that evaluate the numerical fluxes precisely, ensure the stability of the system and keep the physical constraints of the model~\cite{Dedneretal:2002,Kusano:2005}. In particular, we study here the resistive MHD model, which is more realistic in the context of space physics. 
The resistivity adds diffusive effects to the system and allows the simulation of physical events such as magnetic reconnection, a phenomenon which happens, e.g., \textcolor{black}{when the interplanetary magnetic field is merging in particular regions with} the Earth's magnetosphere.

Our goal is to present the verification of the developed framework, that includes a combination of numerical schemes, and the influence, efficiency and stability of the adaptive MR approach for solving the ideal and resistive magneto-hydrodynamic equations for different problems, describing a variety of physical situations. {Our motivation for developing adaptive MHD codes is triggered by space weather applications in science and technology for which there is a social demand.  
High fidelity real-time space weather predictions including the different involved physical phenomena and the computational cost are still challenging~\cite{koskinen2017achievements,morley2019challenges}.}

MHD model simulations in the context of the adaptive MR approach were firstly presented for 1D and 2D Riemann problems \cite{domingues2013extended,Gomesetal:2015}, in a two-dimensional ideal model, followed by the Kelvin--Helmholtz instability \cite{gomes2017ideal}. First results with the three-dimensional implementation were also presented for 1D and 2D Riemann problems in \cite{Gomesetal:2018NSC} and compared with results obtained with the FLASH code. The CARMEN--MHD code for ideal MHD was used in the AMROC framework \cite{lopes2018ideal} and a comparison of the results was performed later in \cite{domingues2019wavelet}, with a wavelet-based adaptive approach. In this work, we present a revised CARMEN--MHD code with new features, which is fully 3D, including resistive terms and which allows the simulation of different problems. Using this new implementation, we present here 2D and 3D resistive and ideal simulations, such as the magnetic reconnection and shock cloud problems. The results are compared with reference solutions.

The remainder of the manuscript is organized as follows.
In Section~\ref{sec1}, we briefly present the resistive MHD model we adopted in its quasi-conservative form and the MHD variables. 
The numerical approach, including the adaptive MR for cell averages, the divergence cleaning and the reference solution information, is presented in Section~\ref{sec2}. In Section~\ref{sec:new}, we describe the implementation and the developed open source code in detail. In Section~\ref{sec3} the numerical results are presented and discussed, and comparisons with reference solutions are given. Conclusions are drawn in Section~\ref{sec4}.

\section{MHD model}
\label{sec1}
We present the MHD model and consider the single fluid description of a plasma, i.e., neglecting the individual identity of each particle of the fluid. 
We are interested in the quasi-conservative form of the MHD model, which expresses local and global conservation of mass and momentum, and quasi-conservation of energy density and magnetic flux. The resistive MHD model in its quasi-conservative form is given by
\begin{subequations}
\begin{eqnarray}
    \frac{\partial \rho}{\partial t} + \nabla\cdot(\rho\textbf{u}) &=&0,\label{eq:MassCons}\\
    \frac{\partial(\rho\textbf{u})}{\partial t} + \nabla\cdot\left[ \rho\textbf{u}^\text{t} \textbf{u} + \left( p +\frac{|\textbf{B}|^2}{2}\right)\textbf{I} -\textbf{B}^\text{t}\textbf{B}\right]&=&0,\label{eq:momCons}\\
    \frac{\partial \mathcal{E}}{\partial t} + \nabla\cdot\left[ (\mathcal{E}+p)\textbf{u} + \textbf{u}\cdot\left( \frac{|\textbf{B}|^2}{2}\textbf{I} - \textbf{B}^\text{t}\textbf{B} \right) \right] &=& \nabla\cdot\left[ \textbf{B}\times\eta(\nabla\times\textbf{B}) \right],\label{eq:enerCons}\\
    \frac{\partial\textbf{B}}{\partial t} + \nabla\cdot\left(\textbf{u}^\text{t}\textbf{B}-\textbf{B}^\text{t}\textbf{u}\right) &=& -\nabla\times(\eta\nabla\times\textbf{B})\label{eq:magCons},\\
    \nabla\cdot\textbf{B} &=&0\label{eq:divB},
\end{eqnarray}
\label{sys:MHDcons}
\end{subequations}
\noindent 
where the magnetic field $\textbf{B}=(B_x,B_y,B_z)$ is the electromagnetic variable, and the fluid variables are the mass density $\rho$, pressure $p$ and velocity $\textbf{u}=(u_x,u_y,u_z)$. 
Without any loss of generality, we define the magnetic field as $\textbf{B}=\textbf{B}/\mu_0$, where $\mu_0$ is the permeability of free space. It is important to note that the MHD variables are normalized. The scalar resistivity is denoted by $\eta=\eta(x,y,z)$, and $\gamma$ is the adiabatic constant. The energy density $\mathcal{E}$ is given by the constitutive law
\begin{equation}
    \mathcal{E} = \frac{p}{\gamma - 1} + \frac{\rho u^2}{2} + \frac{B^2}{2},
\end{equation}
and depends thus on the other variables. The Equations~(\ref{eq:MassCons})~and~(\ref{eq:momCons}) describe the conservation of mass and momentum, respectively. On the other hand, Equations~(\ref{eq:enerCons})~and~(\ref{eq:magCons}) describe the quasi conservation of energy density and magnetic flux. Equation~(\ref{eq:divB}) is the magnetic field constraint, which ensures $\textbf{B}$ is divergence free in the continuous setting. In the absence of resistivity, i.e., $\eta=0$, the model is called the ideal MHD model, which describes the dynamics of a perfectly conducting fluid. In this case, each equation of the system acts as a conservation law, and there is no source term on the right-hand side of the equations. 
%

\section{Numerical approach}
\label{sec2}
In this section we present the numerical methods used in this work. We start with the space discretization of the model, which uses a finite volume formulation. We also recall the multiresolution approach used to adapt the computational mesh and the thresholding of the wavelet coefficients. Equation~(\ref{eq:divB}) of the MHD model is not satisfied numerically, since we are considering a discretized version of the problem. In order to fix it, we use a divergence cleaning, which is also presented in this section.

To introduce the numerical approach of the MHD equations, we first rewrite the System (\ref{sys:MHDcons}) in its vector form
\begin{equation}
    \frac{\partial \textbf U}{\partial t} + \nabla\cdot\mathcal{F}(\textbf U) = \mathcal{S}(\textbf U),
    \label{eq:ConsLaw}
\end{equation}
where $\textbf{U}=(\rho,\rho \textbf{u},\mathcal{E},\textbf{B})$ is the vector of conservative variables, $\mathcal{F} = \mathcal{F}(\textbf U)$ the flux tensor and $\mathcal{S} = \mathcal{S}(\textbf U)$ the vector of source terms, described in Equation~(\ref{sys:MHDcons}). 

\subsection{Finite Volume Discretization}
The Finite Volume (FV) method is based on the integral form of conservation laws, which define the rate of change of a quantity in a fixed volume \cite{Leveque:2002}. To this end the 3D computational domain $\Omega$ is divided into grid cells of the form \[C_{i,j,k}\equiv [x_{i-1/2,j,k}, x_{i+1/2,j,k}]\times[y_{i,j-1/2,k}, y_{i,j+1/2,k}]\times[z_{i,j,k-1/2}, z_{i,j,k+1/2}],\] with $i,j,k\in\{ 0,\cdots,N-1 \}$, where $N$ is the number of cells in each direction and $(x_{i,j,k},y_{i,j,k},z_{i,j,k})$ is the center of the cell $C_{i,j,k}$. In each cell center we define a corresponding cell average, given by
\begin{equation}
\overline{\textbf{U}}_{i,j,k} = \frac{1}{|C_{i,j,k} |}\int_{C_{i,j,k}}{\bf U}(x,y,z,t)\,\text{d}V,
\end{equation}
where $\textbf{U}=\textbf{U}(x,y,z,t)$ is the vector of the variables, $V$ is the volume of the fluid, $|C_{i,j,k}|=\Delta x \Delta y \Delta z$ is the volume of the cell, with $\Delta x = x_{i+1/2,j,k} - x_{i-1/2,j,k}$, $\Delta y = y_{i,j+1/2,k} - y_{i,j-1/2,k}$ and $\Delta z = z_{i,j,k+1/2} - z_{i,j,k-1/2}$. 
By integrating the Equation~(\ref{eq:ConsLaw}) over $C_{i,j,k}$, we obtain
\begin{equation}
\int_{C_{i,j,k}}\frac{\partial{\bf U}}{\partial t}\,\text{d}V + \int_{C_{i,j,k}}\nabla\cdot{\bf F}({\bf U})\,\text{d}V = \int_{C_{i,j,k}}{\bf S}({\bf U})\,\text{d}V.
\label{eq:intcons}
\end{equation}
By multiplying the Equation~(\ref{eq:intcons}) by $\frac{1}{|C_{i,j,k}|}$ and applying the divergence theorem on the divergence operator term, we get
\begin{equation}
\frac{\partial }{\partial t}{ \overline{\bf U }_{i,j,k}} = - \frac{1}{|C_{i,j,k}|}\int_{\partial C_{i,j,k}}{\bf F}({\bf U})\cdot {\bf n}_{i,j,k}\,\text{d}S  + \overline{\bf S}_{i,j,k}({\bf U }),
\label{eq:leiCons}
\end{equation}
where ${\bf n}={\bf n}_{i,j,k}$ is the vector normal to the cell interfaces $C_{i,j,k}$, $\partial C_{i,j,k}$ denotes the boundary of the cell and $dS$ is the surface element of the cell volume. We conclude that the flux tensor must be evaluated on the interfaces of the cell $C_{i,j,k}$, instead of its center. The 2D discrete formulation can be obtained analogously by removing the index $k$.

\subsection{Multiresolution Analysis for Cell-Averages}
The multiresolution representation of cell average data is the essential building block to introduce adaptivity and sparse representation of the solution {\color{black} in the finite volume context}~\cite{Harten:1996}.

To this end, we consider $\textbf{U}=\textbf{U}(x,y,z,t)$ an {\color{black} absolutely}  integrable function on $\Omega$ at a given time instant $t$, and we conceive a hierarchy of dyadic uniform meshes. The center of each cell $C^\ell_{i,j,k}$ is located by $i,j,k\in\{0,\cdots,2^{\ell}-1\}$ and its size is defined as $|C^\ell_{i,j,k}|=h2^{-D\ell}$, where $D$ is the dimension of $\Omega$, $\ell$ the refinement level and $h=|C_{i,j,k}|$. The total number of cells in each level $\ell$ is $2^{D\ell}$ cells. 

Starting with the idea of nested meshes, i.e., a coarser mesh is contained in the finer one, it is possible to navigate between these meshes to obtain cell average values of interest. For this procedure, it is necessary to define two operators: projection and prediction. The projection operator is exact and unique, denoted by $\mathcal{P}_{\ell+1}^{\ell}$, and consists in computing the values on the coarser level $\ell$ from values on the finer level $\ell+1$. This evaluation is accomplished from the weighted average of the cell averages in $\ell +1$, i.e., for three dimensions it is given by
\begin{equation}
\overline{\bf U}^\ell_{i,j,k} = \left(\mathcal{P}_{\ell+1}^{\ell}\overline{\bf U}^{\ell+1}\right)_{i,j,k} = \frac{1}{8}\sum\limits_{m=0}^{1}\sum\limits_{p=0}^{1}\sum\limits_{q=0}^{1}\overline{\bf U}^{\ell+1}_{2i+m,2j+p,2k+q},
\end{equation}
with $m,p,q\in\{0,1\}$. Therefore, each cell average on level $\ell$ is obtained from eight or four values, according to the number of dimensions of the problem.
On the other hand, the procedure that consists in obtaining the cell averages on level $\ell+1$ from the cell averages on the coarser level $\ell$ is performed by the prediction operator, denoted by $\mathcal{P}_\ell^{\ell+1}$. Since this operator predicts the cell average values, it is not exact and can assume different definitions. In this work, we choose the approach proposed by Harten \cite{Harten:1995,harten1994adaptive,HARTI1993153}. Harten's approach {\color{black}has been} extended to two dimensions \cite{bihari1997multiresolution} and three dimensions \cite{Roussel:2003,RSTB03}. The set of approximated cell averages is denoted by $\hat{\textbf{U}}$. Thus, in three dimensions, the approximation is given by
\begin{eqnarray}
\hat{\bf U}^{\ell+1}_{2i+m,2j+p,2k+q} & =& \left( \mathcal{P}_\ell^{\ell+1} \overline{\bf U}^{\ell} \right)_{2i+m,2j+p,2k+q}\nonumber\\  &=& \mathcal{I}( \overline{\bf U}^{\ell},\ell+1,2i+m,2j+p,2k+q),
\end{eqnarray}

where $\mathcal{I}$ is the interpolation operator. 
The prediction operator approximates eight cell averages on each cell in three dimensions, and four in two dimensions. Moreover, it satisfies the localization property, on which the operator only needs the neighbor values to perform the approximation.

For each interpolation we have an associated error, computed by the difference between the cell average $\overline{\textbf{U}}^{\ell+1}$ on level $\ell+1$ and its approximation $\hat{\textbf{U}}^{\ell+1}$, i.e.,
\begin{equation}
{\bf d}^{\ell} = \overline{{\bf U}}^{\ell+1}- \hat{\bf U}^{\ell+1}\quad \rightarrow \quad \overline{{\bf U}}^{\ell+1}=\mathcal{I}(\overline{\bf U}^\ell)+ {\bf d}^{\ell},
\label{eq:error}
\end{equation}
where $\textbf{d}^\ell=\textbf{d}^\ell_{i,j,k}$ is the local error, also called \textit{detail} or \textit{wavelet coefficient}. These coefficients, provide information about the local regularity of the numerical solution. The number of local wavelet coefficients for each approximation varies according to the dimension of the problem: three coefficients for two dimensions and seven for three dimensions.
Thereby, it is possible to establish a one-to-one correspondence $\overline{{\bf U}}^{\ell+1}\longleftrightarrow\{\overline{{\bf U}}^{\ell},{\bf d}^{\ell}\}$. The wavelet coefficients along with the cell averages $\overline{\textbf{U}}^{\ell}_{i,j,k}$, enable us to obtain the cell averages on level $\ell+1$ whenever necessary, resulting in the relation
\begin{equation}
\overline{\textbf{U}}^{\ell+1}_{2i+m,2j+p,2k+q} \longleftrightarrow \{\overline{\textbf{U}}^{\ell}_{i,j,k},\textbf{d}^\ell_1,\textbf{d}^\ell_2,\cdots, \textbf{d}^\ell_{7}\},
\end{equation}
where $m,p,q\in\{0,1\}$ varies in order to obtain every local cell average on level $\ell+1$ in three dimensions. By defining the set of every wavelet coefficient obtained on a local approximation on level $\ell$ as 
    \begin{equation}
    	{\bf D^\ell}=\{ \textbf{d}^\ell_m,\; 1\leq m \leq 2^{D\,\ell}-1 \},
    	\label{detailVector}    	
    \end{equation}   
    the relation can be generalized to the entire mesh, reaching a one-to-one correspondence between $\overline{\bf U}^{\ell+1}$ and $\{{\bf D}^{\ell}, \overline{\bf U}^{\ell} \}$, given by
\begin{equation}
\overline{\bf U}^{L} \longleftrightarrow \{ {\bf D}^{L-1},{\bf D}^{L-2},\cdots,{\bf D}^{0}, \overline{\bf U}^{0} \},
\end{equation}
which characterizes the process of the multiresolution transform operator $\textbf{M}$, defined as
  \begin{equation}
      \overline{\textbf{U}}_{MR}={\bf M}\;\overline{\textbf{U}}^L\;,\qquad\overline{\textbf{U}}^L={\bf M}^{-1}\overline{\textbf{\textbf{U}}}_{MR}.
     \end{equation}
     where $\overline{\textbf{U}}_{MR}=({\bf D}^{L-1},{\bf D}^{L-2},\cdots,{\bf D}^{0}, \overline{\bf U}^{0})$. 
The MR transform satisfies the properties of localization, polynomial cancellation and stability \cite{DGRSESAIM:2011}. The stability guarantees that small perturbations introduced on transformed data on any scale are not uncontrollably amplified in the iterative applications of the multilevel transform. Therefore, we conclude the information of every cell average on the coarser level and the wavelet coefficients of all levels is equivalent to the information of the cell averages on the most refined level.

\begin{figure}[htb!]
\psfrag{L1}{$\ell+1$}
\psfrag{L}{$\ell$}
\psfrag{ij}{\scriptsize{$\overline{\textbf{U}}^{\ell+1}_{i,j}$}}
\psfrag{i1j}{\scriptsize{$\overline{\textbf{U}}^{\ell+1}_{i+1,j}$}}
\psfrag{2i2j}{\scriptsize{$\overline{\textbf{U}}^{\ell+1}_{2i,2j}$}}
\psfrag{2i12j}{\scriptsize{$\overline{\textbf{U}}^{\ell+1}_{2i+1,2j}$}}
\psfrag{2i12j1}{\scriptsize{$\overline{\textbf{U}}^{\ell+1}_{2i+1,2j+1}$}}
\psfrag{2i2j1}{\scriptsize{$\overline{\textbf{U}}^{\ell+1}_{2i,2j+1}$}}
\centering
\includegraphics[width=0.7\textwidth]{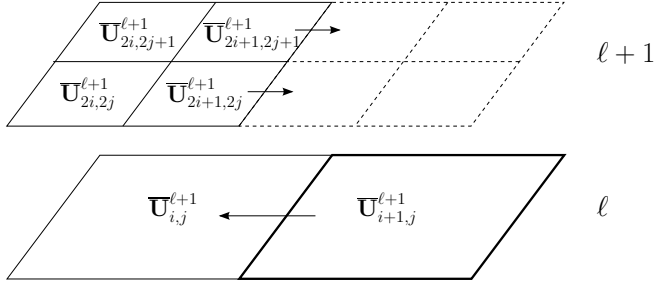}
\caption{\label{fig:fluxes}
{Conservative flux computation in 2D for two different refinement levels, illustrating ingoing and outgoing fluxes. Adapted from \cite{RSTB03}.}}
\end{figure}

{
The conservation properties of the finite volume method are preserved in our adaptive multiresolution discretization taking special care in the flux evaluation.
To ensure the balance of ingoing and outgoing fluxes at the cell interfaces on adjacent refinement levels, we use the conservative formulation proposed in~\cite{RSTB03,Mueller:2003}. To guarantee the conservation the ingoing fluxes at level $\ell$ are computed as the outgoing fluxes of the corresponding cells at level $\ell + 1$, as illustrated in Figure~\ref{fig:fluxes}. This is possible due to the graded-tree structure, which keeps the nearest cousins of a cell or creates a virtual leaf for the flux computations. Thus the flux computation is conservative between cells at different levels of refinement.}

\section*{Dynamic Mesh Adaptation: Thresholding}
To define the regions of the computational mesh that need more or less refinement, we apply the threshold operator $\mathcal{T}_{\epsilon^\ell}$ on the wavelet coefficients. For each approximation, this nonlinear operator is defined as
\begin{equation}
  	\mathcal{T}_{\epsilon^\ell}({\bf d^\ell_m})=
  	\begin{cases}
  	0,\; \text{if}\;|{\bf d}^\ell_m|\leq\epsilon^\ell,\\
  	{\bf d^\ell_m},\;\text{otherwise},
  	\end{cases}
	\label{threshold}
\end{equation} 
where $\epsilon^\ell$ is the threshold parameter and $m\in\{1,\cdots,2^{D\,\ell}-1\}$. Given the level $\ell$, the thresholding consists in removing the cells in which the magnitude of the details is smaller than $\epsilon^\ell$, replacing them by zero. Thereby, the number of cells required for the numerical simulation can be significantly decreased, impacting the computational cost which correspondingly decreases with this reduction. Hence the computational cost becomes smaller when more wavelet coefficients are removed.
In regions where the solution presents smooth behavior, the wavelet coefficients have small magnitude ($|{\bf d}^\ell|\leq\epsilon^\ell$), allowing locally coarser meshes. On the other hand, the magnitude of the coefficients is significant ($|{\bf d}^\ell|>\epsilon^\ell$) in regions where local structures are present, requiring more refined meshes \cite{HARTI1993153,harten1994adaptive}.
 
The threshold parameter $\epsilon^\ell$ can have {\color{black} either} a fixed value or be level dependent. In the former case, we define the value $\epsilon=\epsilon^\ell$ which remains the same during the simulation. In the level dependent case, an initial parameter $\epsilon^0$ is defined and it changes according to the local refinement of the region of interest, given by the equation
\begin{equation}
\epsilon^\ell=\frac{\epsilon^0}{|\Omega|} 2^{D(\ell-L+1)}, \;\;0\leq \ell \leq L-1,
\label{eq:harten}
\end{equation}
where $|\Omega|$ is the global volume of the computational region and $L$ is the maximum refinement level. Equation~(\ref{eq:harten}) is called Harten's strategy to determine the choice of {\color{black} the threshold} parameter $\epsilon^\ell~=~\epsilon(\epsilon^0,\ell)$.

Each detail $\textbf{d}^\ell$ is understood as a vector with the details of the MHD conservative variables as its components, i.e., $\textbf{d}^\ell~=~(\textbf{d}^\ell_\rho,\textbf{d}^\ell_{\rho \textbf{u}},\textbf{d}^\ell_{\mathcal{E}},\textbf{d}^\ell_{\textbf{B}})$. As the MHD variables are stored in the vector $\overline{\bf U}^\ell_{i,j,k}$, the detail components are computed from the approximation error in each variable and its maximum value.
In the scalar-valued approach, the details are computed as $\textbf{d}^\ell_{\varpi}=\frac{\textbf{d}^\ell_{\varpi}}{\max(|\varpi|)}$, where $\varpi$ denotes a MHD quantity. 
In this case, each component of the vector variables $\rho\textbf{u}$ and $\textbf{B}$ is considered separately to compute $\textbf{d}^\ell_{\varpi}$, as scalar variables. Thus the vector of details $\textbf{d}^\ell$ has in total 8 components. 

After performing some {\color{black} numerical} experiments, we found the optimal computation of the details of the MHD variables in three dimensions, given by the following,
  \begin{eqnarray}
  				\textbf{d}^\ell_{\rho}&=&\frac{\textbf{d}^\ell_{\rho}}{\max|\rho|},\vspace{5mm}\\ 
  				\textbf{d}^\ell_{\rho\textbf{u}}&=&\frac{\sqrt{(\textbf{d}^\ell_{\rho u_x})^2 + (\textbf{d}^\ell_{\rho u_y})^2 + (\textbf{d}^\ell_{\rho u_z})^2}}{\max(|\rho u_x|,|\rho u_y|,|\rho u_z|)},\vspace{5mm}\\
  				\textbf{d}^\ell_{\mathcal{E}}&=&\frac{\textbf{d}^\ell_{\mathcal{E}}}{\max|\mathcal{E}|},\vspace{5mm}\\
  				\textbf{d}^\ell_{\textbf{B}}&=&\frac{\sqrt{(\textbf{d}^\ell_{B_x})^2 + (\textbf{d}^\ell_{B_y})^2 + (\textbf{d}^\ell_{B_z})^2}}{\max(|B_x|,|B_y|,|B_z|)}. 
  		\end{eqnarray}		
In two dimensions we have a special treatment for the $z$-components, namely 
  		\begin{eqnarray}
  				\textbf{d}^\ell_{\rho\textbf{u}}&=&\frac{\sqrt{(\textbf{d}^\ell_{\rho u_x})^2 + (\textbf{d}^\ell_{\rho u_y})^2 }}{\max(|\rho u_x|,|\rho u_y|)},\vspace{5mm}\\
  				\textbf{d}^\ell_{\rho u_z}&=&\frac{\textbf{d}^\ell_{\rho u_z}}{\max|\rho u_z|},\vspace{5mm}\\ 
  				\textbf{d}^\ell_{\textbf{B}}&=&\frac{\sqrt{(\textbf{d}^\ell_{B_x})^2 + (\textbf{d}^\ell_{B_y})^2}}{\max(|B_x|,|B_y|)},\vspace{5mm}\\
  				\textbf{d}^\ell_{B_z}&=&\frac{\textbf{d}^\ell_{B_z}}{\max|B_z|},\vspace{5mm}
  		\end{eqnarray}
which is called vector-based approach. Therefore, for vector-valued variables $\rho\textbf{u}$ {\color{black} and} $\textbf{B}$, the associated wavelet coefficients take into account each component on the computation. In two dimensions, only two components are used and the third one is computed in a independent way. We found that this vector-based approach optimizes the local mesh refinement, compared to the scalar-valued approach, where we compute the details for each vector component individually~\cite{GomesThesis2017}.

{The adaptive mesh is organized into a tree data structure where the different levels define the resolution, which are represented in the tree hierarchically. 
Following the tree nomenclature, the children of a given cell $C=C^\ell$ are the cells descendent of $C$ in level $\ell -1$. 
The neighborhood of a cell at the same level corresponds to their brothers and
in the upper level to their uncles.
The adaptive mesh must follow a graded structure as discussed in \cite{Cohen:2003book}. This contributes to maintain the accuracy and stability in the time evolution with the inclusion of extra neighbors and uncles cells in the neighborhood of the selected cells. We also add virtual leaves to avoid unnecessary computation. Virtual leaves are not evolved in time.  
More details can be found in~\cite{RSTB03}.
The following Algorithm~\ref{alg-mesh} describes the procedure to construct the adaptive mesh at $t=t_0$. The adaptive mesh is then updated at each time step.}

\begin{figure}[H]
\centering
\caption{2D MHD shock-cloud problem: adaptive meshes obtained at the final time of simulation with (a) scalar-value and (b) vector-based threshold approaches.} 
\begin{tabular}{cc}
 \hspace{6mm}(a) &\hspace{6mm}(b)\\
\hspace{8mm}\includegraphics[width=0.3\textwidth]{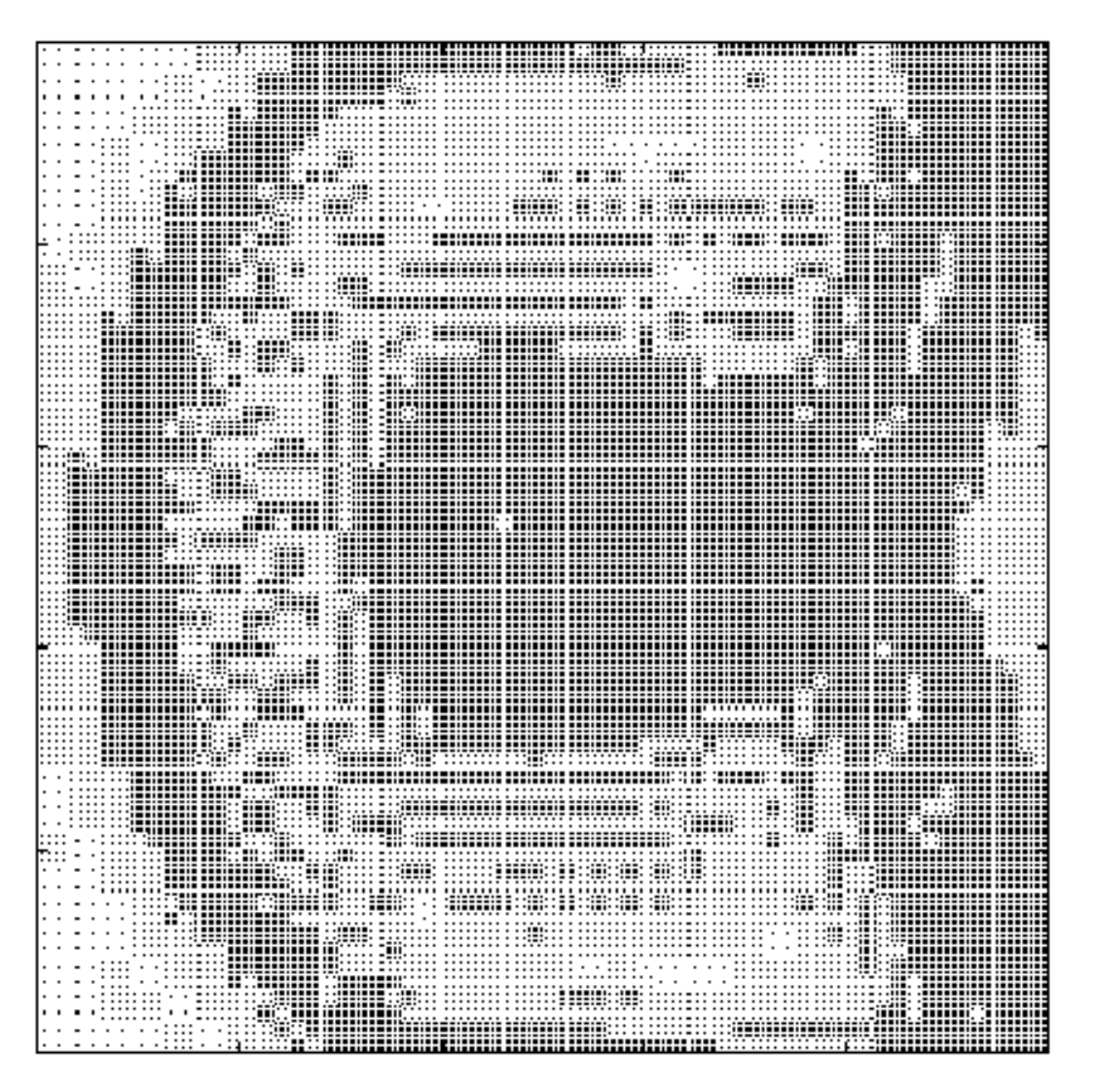}&\hspace{8mm}\includegraphics[width=0.3\textwidth]{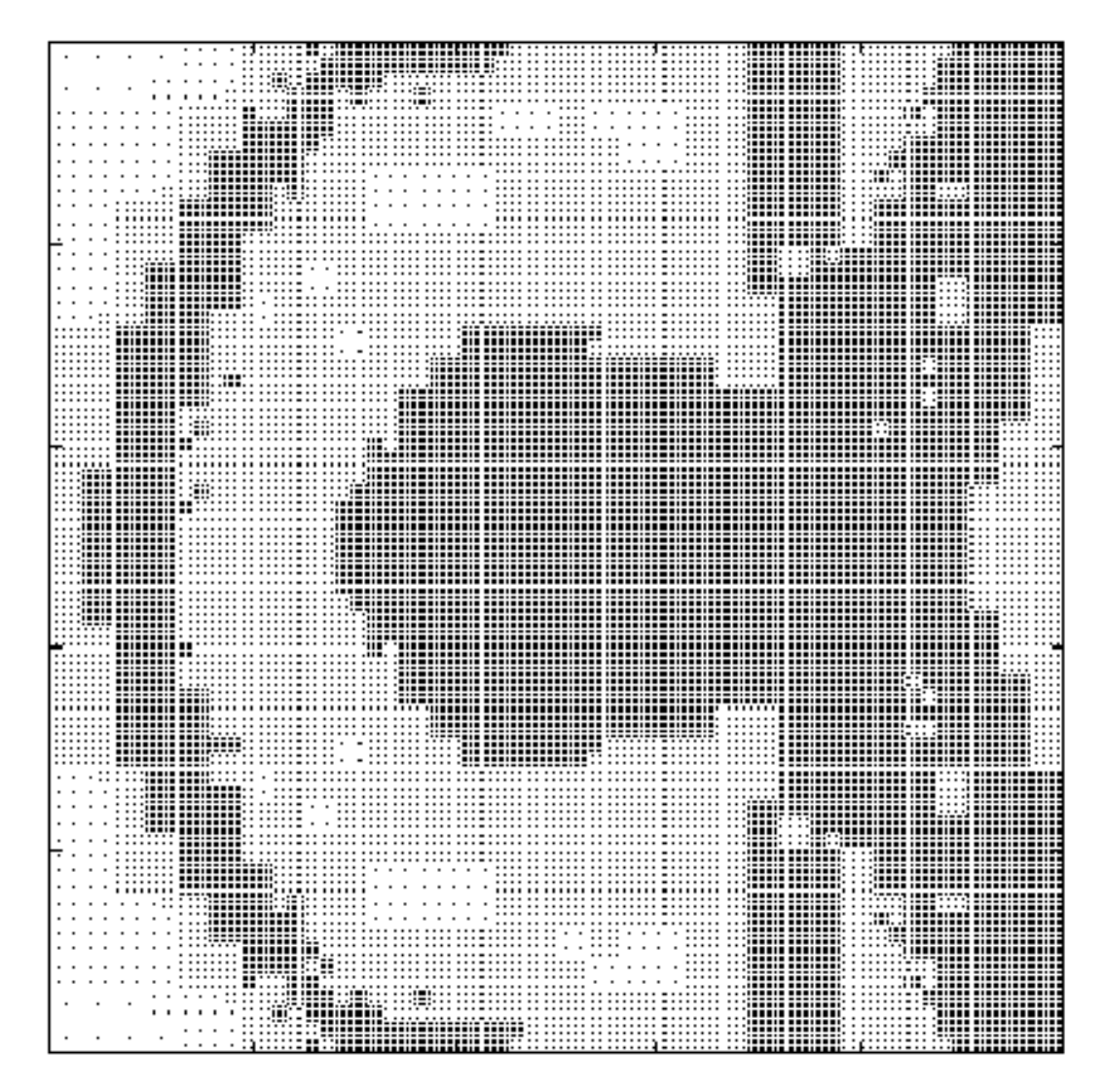}
\end{tabular}
\label{fig:meshes}
\end{figure}

\begin{algorithm}
\caption{Adaptive Mesh Construction}
\begin{algorithmic}
\label{alg-mesh}
\REQUIRE Consider the adaptive mesh, containing all cells at level $L$
\REQUIRE Threshold parameter $\epsilon^\ell$

\STATE Obtain the set of details $\textbf{D}^\ell$ of every level $\ell$
\FOR{$\ell = L-1;\,\, \ell \ge \ell_{\text{min}};\,\, \ell = \ell-1$}
\FOR{Every cell $\mathcal{C}$ at level $\ell$}
\IF{Every child of $\mathcal{C}$ is in the adaptive mesh}
\IF{The significant details are kept, \textit{i.e.} the detail correspondents to cell $\mathcal{C}$ is smaller then $\epsilon$}
\IF{The refinement level of every cell adjacent to the set of cells to be merged is greater than $\ell$}
\STATE Remove the child of $\mathcal{C}$ from the adaptive mesh
\STATE Insert $\mathcal{C}$ in the adaptive mesh
\STATE Check if the tree is graded, if not impose it.
\STATE Add virtual leaves.
\ENDIF
\ENDIF
\ENDIF
\ENDFOR
\ENDFOR 
\STATE Refine every cell at a level $\ell<L$;
\end{algorithmic}
\end{algorithm}
In Figure~\ref{fig:meshes}, we present two adaptive meshes obtained at the final time for 2D shock-cloud simulations for the same $\epsilon$ using either the (a) scalar-valued or (b) a vector-based threshold approach. By comparing the adaptive meshes with the solutions of the variables presented in Figure~\ref{fig:SC2D1}, computed with the vector-based approach on a adaptive mesh, we can observe that the vector-based approach is more efficient to capture the structures of the numerical solution, avoiding unnecessary refinement in smoother regions of the solution. 

To apply the threshold operator to previously defined values, we define the maximum value among the details, i.e., $\max|\textbf{d}^\ell| = \max\{|\textbf{d}^\ell_\rho|,|\textbf{d}^\ell_{\rho \textbf{u}}|,|\textbf{d}^\ell_{E}|,|\textbf{d}^\ell_{\textbf{B}}|\}$. 
Thus, after this procedure, the threshold operator is applied to the details of the MHD variables. 
In this way, the adaptive mesh becomes the union of {\color{black} the mesh of each variable}, since the wavelet coefficients of each variable are used to decide in which local region a more refined mesh is necessary.

\subsection{Divergence of the Magnetic Field Correction}

Gauss' law of magnetism, given by Equation~(\ref{eq:divB}), imposes a physical constraint on the magnetic field, which is satisfied in the 
{\color{black} continuous}
medium. To ensure the absence of nonphysical behavior in the numerical MHD solution, we add to the MHD system the so-called parabolic-hyperbolic divergence cleaning \cite{Dedneretal:2002,mignone2010second}, which does not impose a vanishing divergence of the magnetic field, but {\color{black} damps} and propagates the associated numerical divergence errors. In this case, we add a new scalar variable $\psi$ to the MHD System~(\ref{sys:MHDcons}), thus modifying Equation~(\ref{eq:magCons}) and adding a transport equation for $\psi$, 
\begin{eqnarray}    
		\frac{\partial{\textbf B}}{\partial t} +\nabla\cdot\left( \textbf{u}^\text{t}\textbf{B} - \textbf{B}^\text{t}\textbf{u} \right) + \nabla\psi &=& -\nabla\times\left({\eta\nabla\times{\textbf B}}\right), \label{eq:Mag}\\
		\frac{\partial \psi}{\partial t} + c_h^2\nabla\cdot{\textbf B} &=&-\frac{c_h^2}{c_p^2}\psi,
		\label{eq:divCorr}
\end{eqnarray}
where $c_p$ and $c_h$ are the parabolic and hyperbolic constants, with $c_h>0$. 
The complete model is called MHD model with Generalized Lagrange Multipliers (GLM--MHD). It is important to note that the GLM--MHD model is originally proposed for ideal MHD and we use this formulation here too. In the resistive case, we are only considering the additional source terms present Equations~(\ref{eq:enerCons})~and~(\ref{eq:magCons}). 

By using appropriated initial and boundary conditions, the MHD system presented is completed and ready for the numerical simulation in two or three space dimensions.

\section{Time Evolution}
{
The adaptive mesh has to be updated at each time step, because the local structures present in the numerical solution can change at each iteration.
By defining the operators $\mathbb{E}:\overline{\textbf{U}}^n\rightarrow\overline{\textbf{U}}^{n+1}$ of the time evolution, $\mathcal{T}: \{ \textbf{d}^{L-1},\cdots,\textbf{d}^\ell,\overline{\textbf{U}}^\ell \}\to \overline{\textbf{U}}^{n\star}$ the thresholding operator, $\textbf{M}:\overline{\textbf{U}}^{L}\rightarrow\{ \textbf{d}^{L-1},\cdots,\textbf{d}^\ell,\overline{\textbf{U}}^\ell \}$ of the multiresolution transform and $\textbf{M}^{-1}:\{ \textbf{d}^{L-1},\cdots,\textbf{d}^\ell,\overline{\textbf{U}}^\ell \}\rightarrow \overline{\textbf{U}}^{L}$ of the inverse multiresolution transform, the adaptive MR process can be described as follows
\begin{eqnarray}
\overline{\bf U}^{n}_{MR}   &=& {\bf M}(\overline{\bf U}^{n})\\
\overline{\bf U}^{n\star}_{MR} &=& \mathcal{T}(\overline{\bf U}^n_{MR}),\\
\overline{\bf U}^{n+1}_{MR} &=& \mathbb{E} (\overline{\bf U}^{n\star}_{MR}),\\
\overline{\bf U}^{n+1}      &=& {\bf M}^{-1}(\overline{\bf U}^{n+1}_{MR}).
\end{eqnarray}
This process indicates that, after the MR representation of the solution $\overline{\bf U}^{n}$, we apply the thresholding operator $\mathcal{T}$ and obtain the cell averages at the intermediate step $n^\star$. Then, these values are evolved to the time $t^{n+1}$. To finish the process, the inverse MR operator is performed~\cite{RSTB03}.}

\section{Implementation issues and the CARMEN--MHD code}
\label{sec:new}
The CARMEN--MHD code is based on the CARMEN code, originally developed by O. Roussel during his PhD thesis~\cite{Roussel:2003,RSTB03} using finite volumes together with adaptive MR for cell-averages. This $C^{++}$ code with tree-data structures was implemented to simulate reaction-diffusion equations modeling combustion problems and later also extended for the compressible Euler and Navier--Stokes equations. 

The implementation of the CARMEN--MHD code started with the ideal 2D MHD equations, with HLL and HLLD numerical fluxes and GLM divergence cleaning \cite{Gomes:2012:AnMuAd,domingues2013extended}. After some adjustments, including an eigenvalue fix, TVD limiters for the conservative variables and improvement of the boundary conditions  \cite{Gomesetal:2015,gomes2017ideal}, we started the implementation of the 3D MHD equations, first for 2.5D simulations \cite{Gomesetal:2018NSC}. The uniform mesh MHD implementation in 2D and 3D, which allowed the CPU time comparison, resistive terms, artificial diffusion terms, fixed time steps and more, were implemented later \cite{GomesThesis2017}. The code became more robust and different types of MHD simulations could be done properly.

The adaptive MR algorithm creates a computational mesh which becomes more refined in regions where local structures are present. The mesh refinement reduces the cost of the numerical flux computation significantly, which normally requires the major memory percentage. The numerical simulation of the ideal and resistive MHD equations is performed with the CARMEN--MHD code. 

To compute the flux $\mathcal{F}$ on the cell interface, we firstly reconstruct the conservative variables using a MUSCL-type monotonized central scheme \cite{van1974towards} to achieve second order accuracy in space. 
These reconstructed values are used to compute the intermediate states of the Harten-Lax-Van Leer-Discontinuities Riemann solver \cite{Kusano:2005} and then the numerical flux is evaluated on the cell interfaces. We should note that the FV method is strictly conservative, since the outflux of a volume is imposed to be equal to its influx. The physics of the problem can be reproduced in such way that the conservative principles of the model are sustained.  

The current CARMEN--MHD code is able to solve ideal and resistive MHD equations, by using the adaptive MR for cell averages or a uniform finite volume computational mesh. The simulations performed are stable and the numerical solution remains coherent and does not show oscillations or non physical behavior. The code and its documentation are available on a repository for download and all the presented problems can be reproduced and visualized properly\footnote{waveletapplications.github.io/carmenMHD/}.

\begin{figure}[H]
\centering
\includegraphics[width=0.95\textwidth]{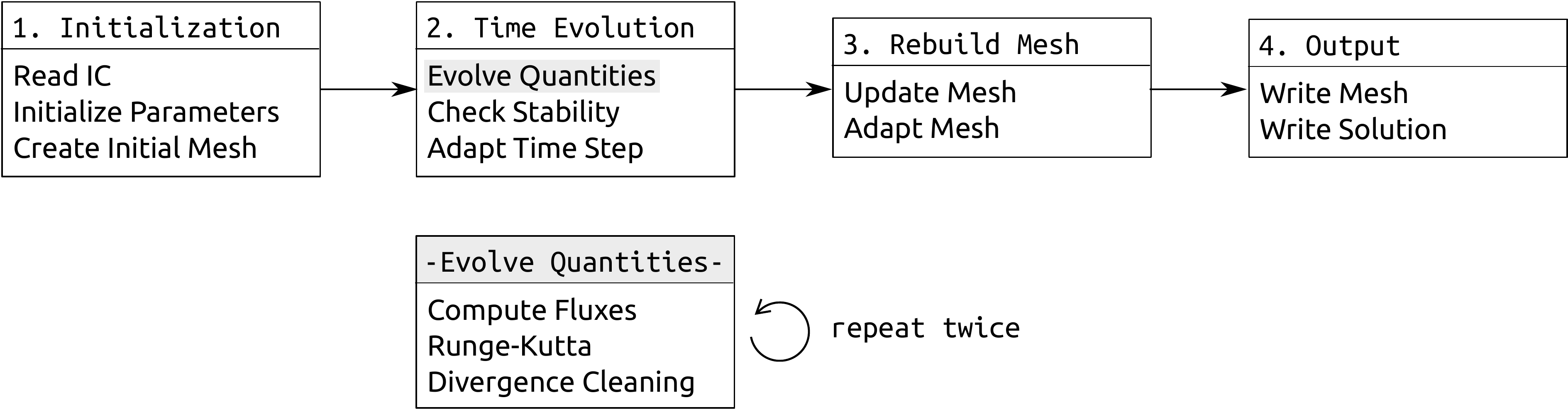}
\caption{{\color{black} Flowchart of the CARMEN--MHD code.}}
\label{fig:fluxogram}
\end{figure}

In Figure~\ref{fig:fluxogram}, we present a flowchart illustrating the algorithm of the CARMEN--MHD code. First, the code is initialized with the initial condition (IC) and other parameters. The cell averages are computed and the initial mesh is created. In the second step, the time evolution is performed. In this part, the code evolves the quantities of the MHD model, by evaluating the numerical fluxes, GLM divergence cleaning and evolving in time with a second order Runge--Kutta scheme. The stability of the new solution is checked and the time step of the next iteration is computed. After that, we proceed to the third step, where the adaptive mesh is updated based on the new cell averages, and the mesh is adapted again. To finish the algorithm, the mesh and solution are written into a file. The procedure is repeated until the final time is reached.

As a reference solution for our results, we use the FLASH code\footnote{flash.uchicago.edu/site/flashcode/}, developed in the \textit{Flash Center} at the University of Chicago. This code includes the implementation of the ideal and resistive (constant resistivity) MHD models, as well as the FV discretization. It is possible to perform adaptive simulations by using an adaptive mesh refinement algorithm. However, here we are interested only in the results obtained on a uniform mesh, as the comparison between two adaptive methods is not part of this work. The FLASH code results are used only for the comparison of the solutions and error computations. The version of the code used here is $4.3$. For the FLASH code simulations, the following settings are used: one-step Hancock for time evolution, \textit{8-wave} divergence cleaning, MC limiter and HLLD Riemann solver. These settings yield second order of the numerical scheme in time and space.

\section{Numerical Simulation}
\label{sec3}
In the following we present several test cases to verify the CARMEN--MHD code, and, in particular, the adaptive multiresolution algorithm for the MHD equations. 
These test cases can assess how our solvers deal with different physical situations, such as, magnetohydrodynamic shocks, local structures and magnetic field lines topology changes, and also numerical challenges, such as stability, strong discontinuities and divergence free correction.

In previous works, we found that the MR approach is efficient to represent the numerical solution 2D and 2.5D for ideal MHD problems~\cite{domingues2013extended,Gomesetal:2015,Gomesetal:2018NSC}. The MR algorithm decreases significantly the number of the cells in the computational mesh and, consequently, the required CPU time. It also provides an accurate solution, compared to the regular mesh solution, demanding much less cells and memory. To verify the solutions of the CARMEN--MHD code, we use the FLASH code \cite{Flash:2000} as reference solution for our results. The reference solutions obtained with FLASH code are simulated using {\color{black} finite volumes on a regular Cartesian mesh}.

\subsection{Orszag--Tang Vortex}
\begin{figure}[H]
\psfrag{PRE}{$p$}
\psfrag{p}{$p$}
\psfrag{x}{$x$}
\caption{Cuts in variable $p$ on a uniform mesh: (a) comparison between reference solution (solid line), Londrillo~e~Del~Zanna 2000 (cross) e Miyoshi~e~Kusano~2005 (dot); (b) CARMEN--MHD code solution obtained for $L=7$ (cross), $L=8$ (dot) and $L=9$ (solid line).}
\centering
\begin{tabular}{cc}
$(a)$ & $(b)$\\
\includegraphics[width=0.45\textwidth]{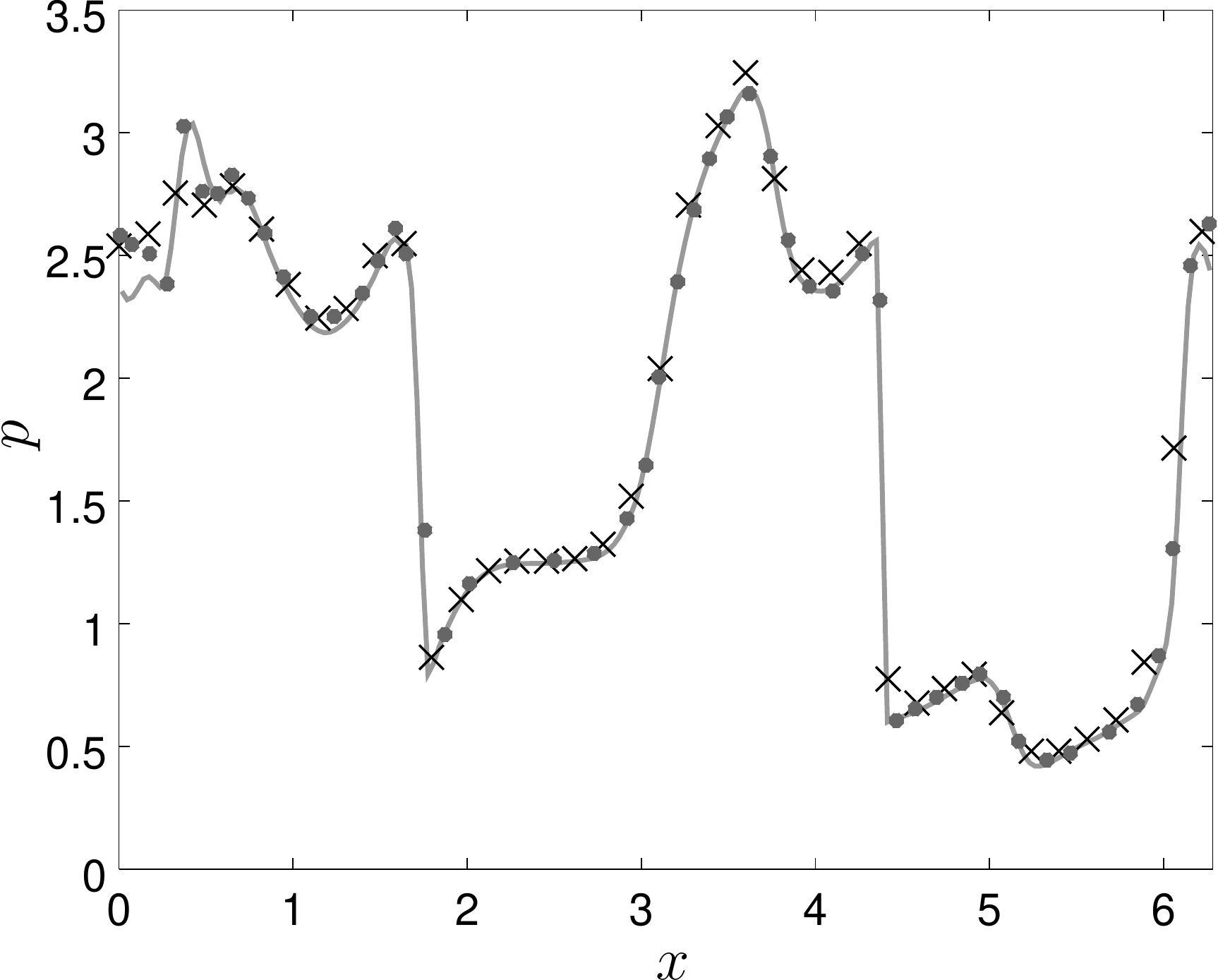}&\includegraphics[width=0.45\textwidth]{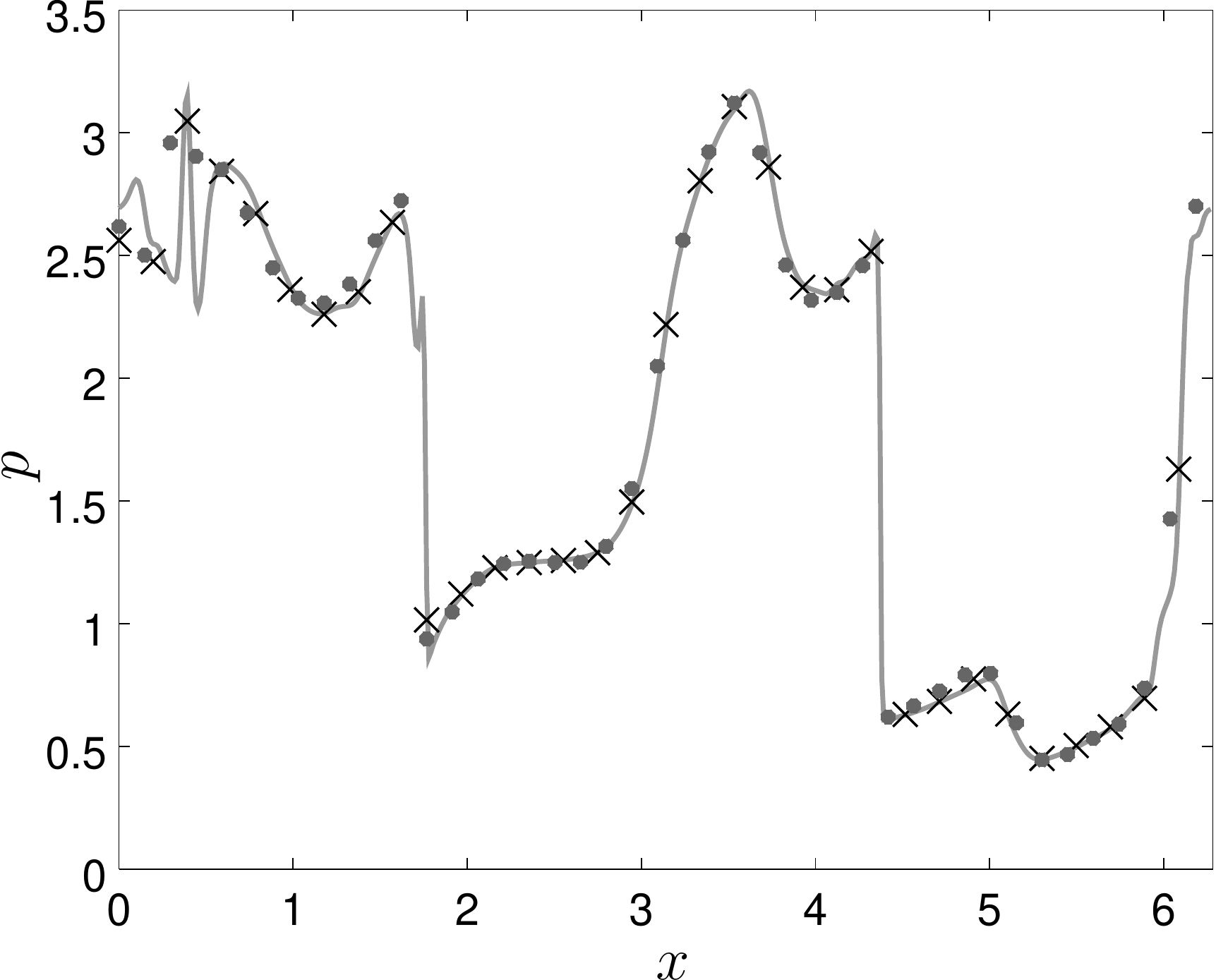}
\end{tabular}
\label{fig:compareCuts1}
\end{figure}

The Orszag--Tang vortex (O-T) in two dimensions~\cite{Orszag1979} is a well-known benchmark for MHD simulations, which allows us to test the transition to two-dimensional supersonic MHD turbulence. Thus, the Orszag--Tang vortex is adequate to verify the robustness of the code when it comes to deal with the formation of magneto-hydrodynamic shocks and shock-shock interactions. It is also interesting to quantitatively estimate how significant the magnetic monopoles affect the numerical solution, by testing the divergence constraint of the magnetic field. In summary, this problem is a common and classical numerical test for MHD codes and consistent to perform comparisons between codes. 
%
This problem presents physical structures over the entire domain, characterizing a challenge to our proposed adaptive multiresolution algorithm. We want to measure the quality of the CARMEN--MHD solution and compare it to the reference. The initial condition for the O-T problem is given in Table~\ref{tab:ot2d}. The domain interval is defined by $\Omega=[0,2\pi]\times[0,2\pi]$, $L=9\,(512\times 512)$ is the most refined level, the final time is $t=\pi$, the Courant number $\nu=0.4$, and $\alpha_p=0.4$ and $\gamma=5/3$. The boundary conditions are periodic {\color{black} in all directions}.
\begin{table}[H]
\caption{Initial condition of the Orszag--Tang vortex. }
\centering
\begin{tabular}{cccccccc}
\toprule
$\rho$ & $p$ & $u_x$ & $u_y$ & $u_z$ & $B_x$ & $B_y$ & $B_z$ \\
\cmidrule{1-8} \\[-2mm]
$\gamma ^2$ & $\gamma$ & $-\sin y$ & $\sin x$ & 0.0 & $-\sin y$ & $\sin 2x$ & 0.0 \\ [2mm]
\bottomrule 
\end{tabular}
\label{tab:ot2d}
\end{table}
Firstly, we present a comparative study of the numerical solution obtained on a regular full mesh, to ensure the reference solution is adequate and to show that CARMEN--MHD results are coherent. 
To evaluate the local convergence of the solution, we present cuts of the variable $p$ at $y=0.64\pi$, which is largely used on the literature, see e.g., \cite{Kusano:2005,londrillo2000high,jiang1999high,ryu1998divergence}. By collecting a set of $200$ points of the solution presented in Londrillo~e~Del~Zanna \cite{londrillo2000high}, Miyoshi~e~Kusano~\cite{Kusano:2005}, along with the FLASH code solution, we can observe in Figure~\ref{fig:compareCuts1}~$(a)$ that these cuts present similar behavior. It also suggests that the reference solution obtained with the FLASH code is adequate to be used as a benchmark for the CARMEN--MHD code solution. The solutions at levels $L=7, 8$ and $9$ obtained with CARMEN--MHD are shown in Figure~\ref{fig:compareCuts1}~$(b)$. For each resolution level, the solution behaves as expected. 

\begin{figure}[H]
\caption{Variables $p$ and $u_y$ at $t=\pi$ and $L=9$, obtained with resistive CARMEN--MHD code for the 2D Orszag--Tang vortex and $\epsilon=0.03$.}
\begin{center}
\begin{tabular}{cc}
$p$ & $u_y$\\
\includegraphics[width=0.44\textwidth]{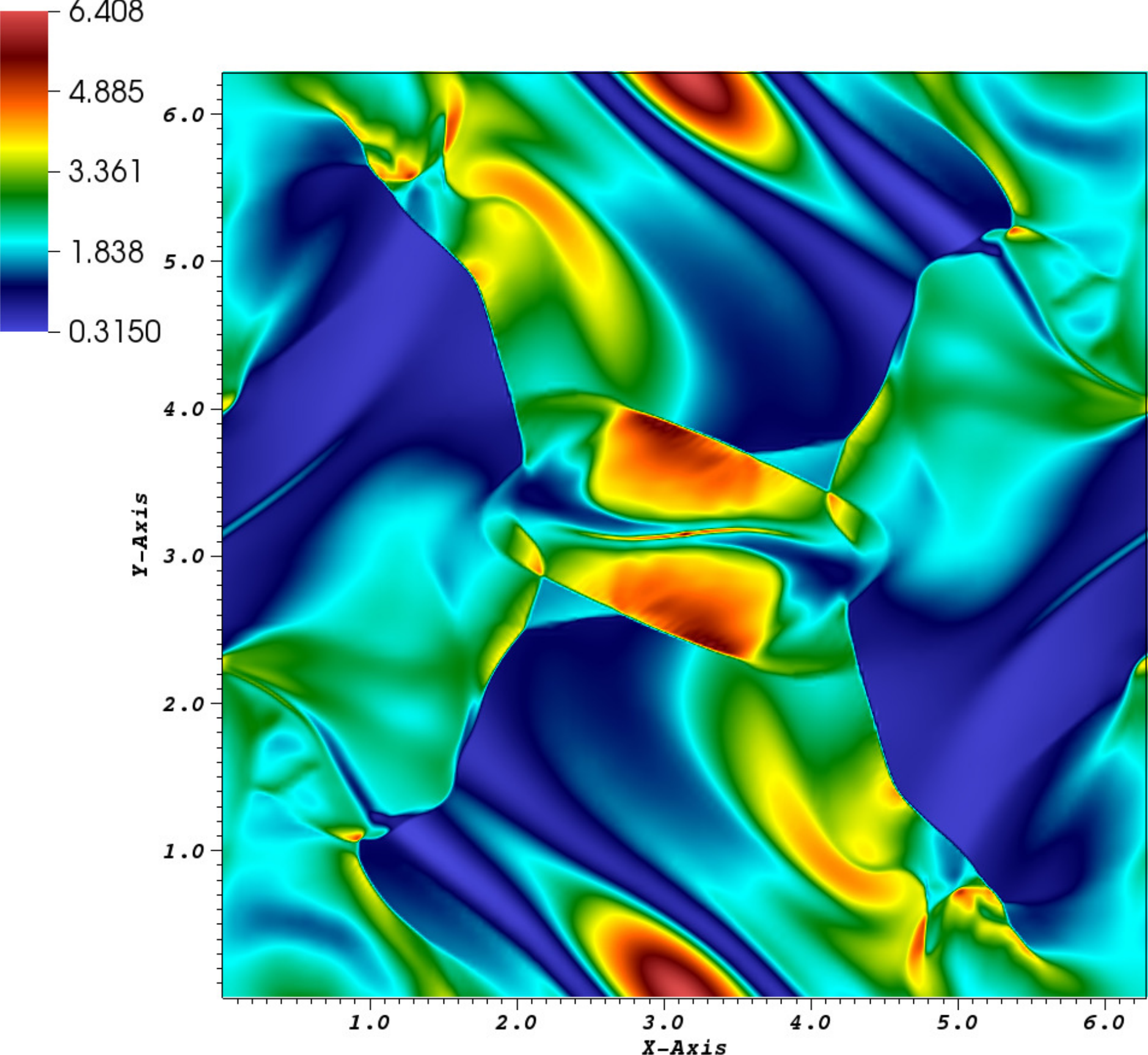} & \includegraphics[width=0.44\textwidth]{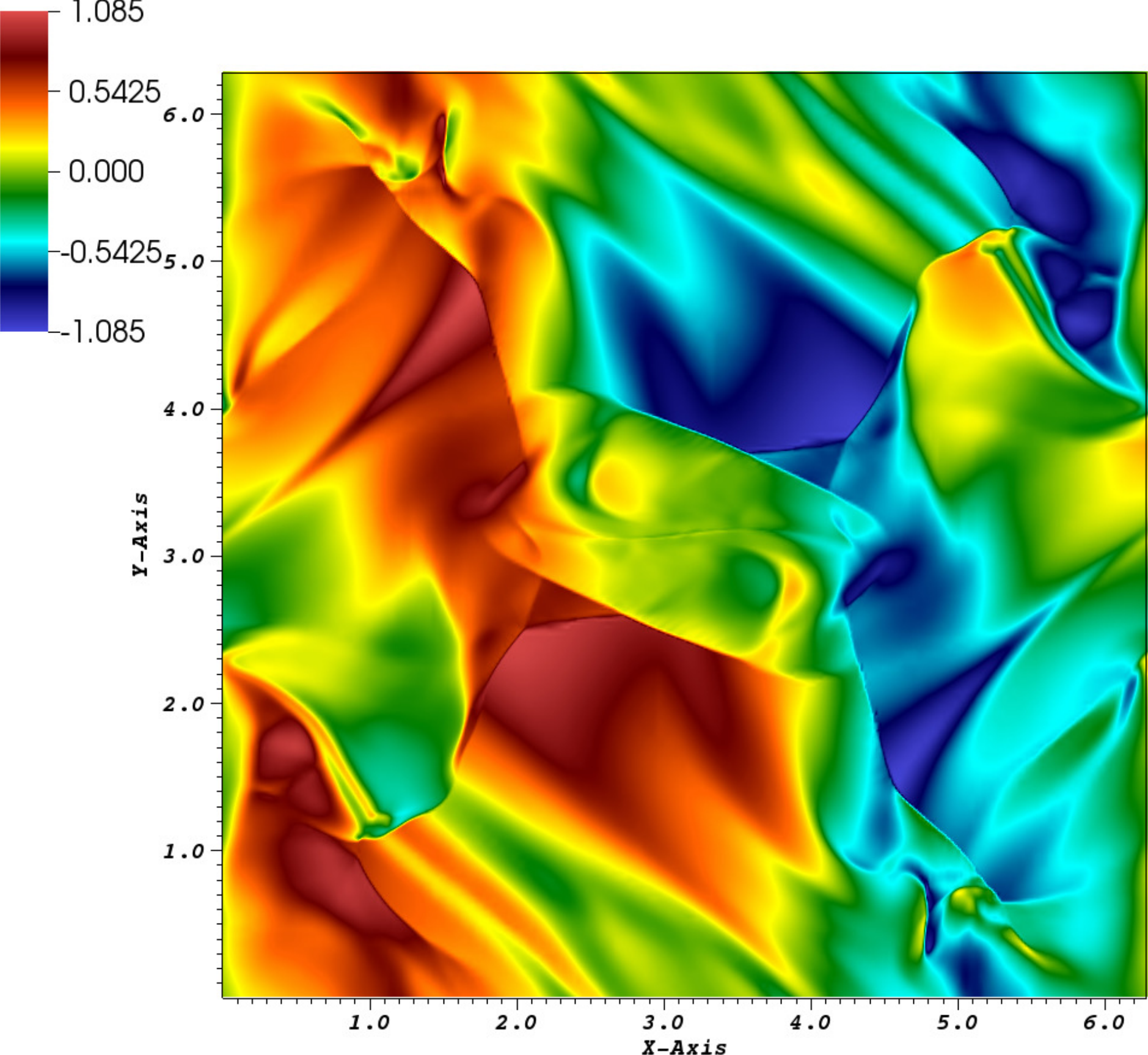} \\
\end{tabular}
\end{center}
\label{fig:OT2DL93}
\end{figure}

\begin{figure}[H]
\psfrag{PRE}{$p$}
\psfrag{VX}{$u_x$}
\psfrag{VY}{$u_y$}
\psfrag{y}{$y$}
\psfrag{x}{$x$}
\psfrag{EPS001}{\scalebox{.5}{$\epsilon=0.01$}}
\psfrag{EPS003}{\scalebox{.5}{$\epsilon=0.03$}}
\psfrag{EPS-02}{\scalebox{.5}{$\epsilon^0=0.2$}}
\caption{Cuts of the variable $p$ (from top to bottom) at $t=\pi$ and $L=9$, for reference solution (line) and CARMEN--MHD adaptive solutions with $\epsilon = 0.01$ (cross), $\epsilon = 0.03$ (circle), $\epsilon^0 = 0.2$ (dot), at  $x=\pi$ and $y=\pi$}
\begin{center}
\begin{tabular}{cc}
$x=\pi$ & $y=\pi$ \\
\includegraphics[width=0.45\linewidth]{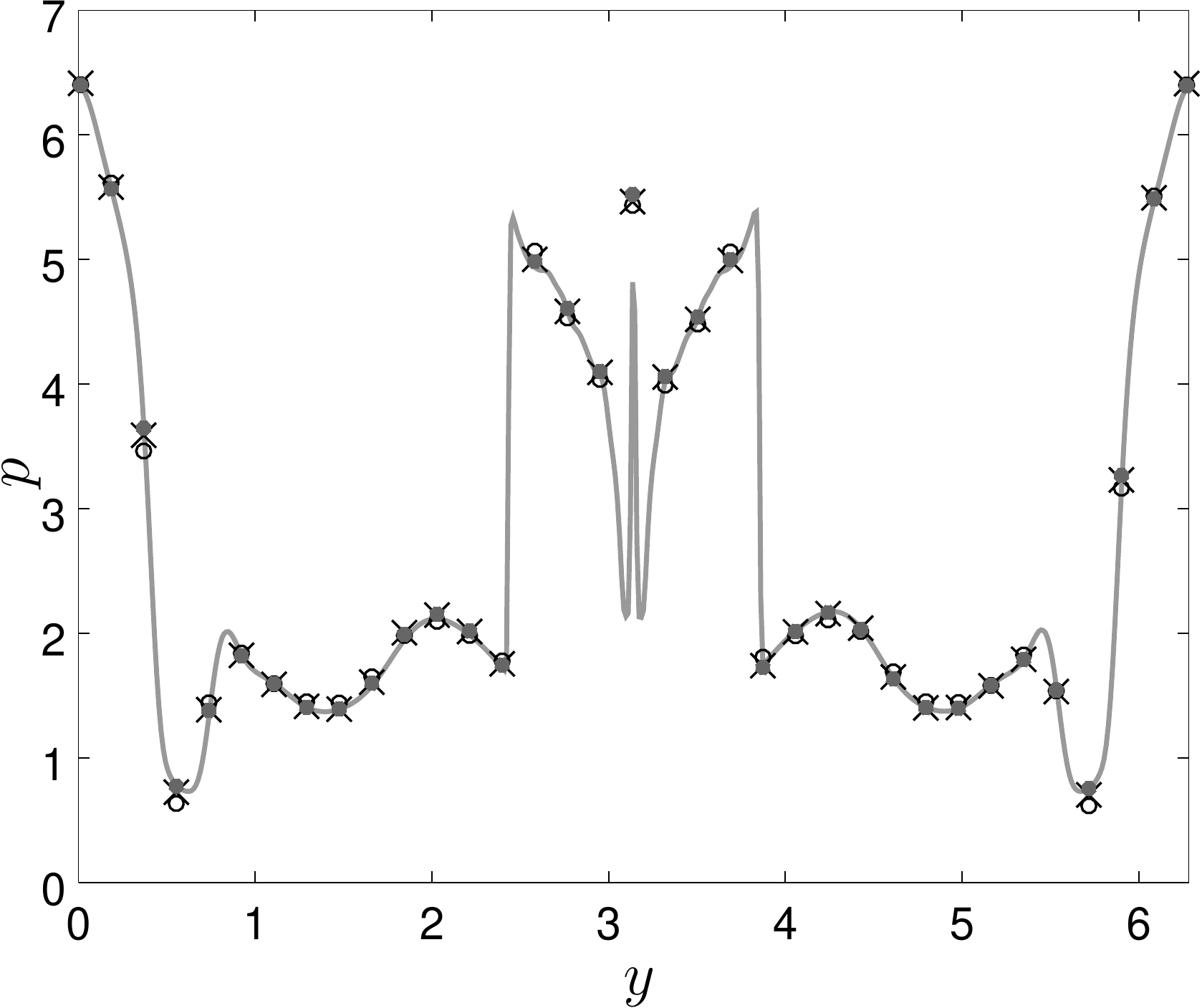} & \includegraphics[width=0.45\linewidth]{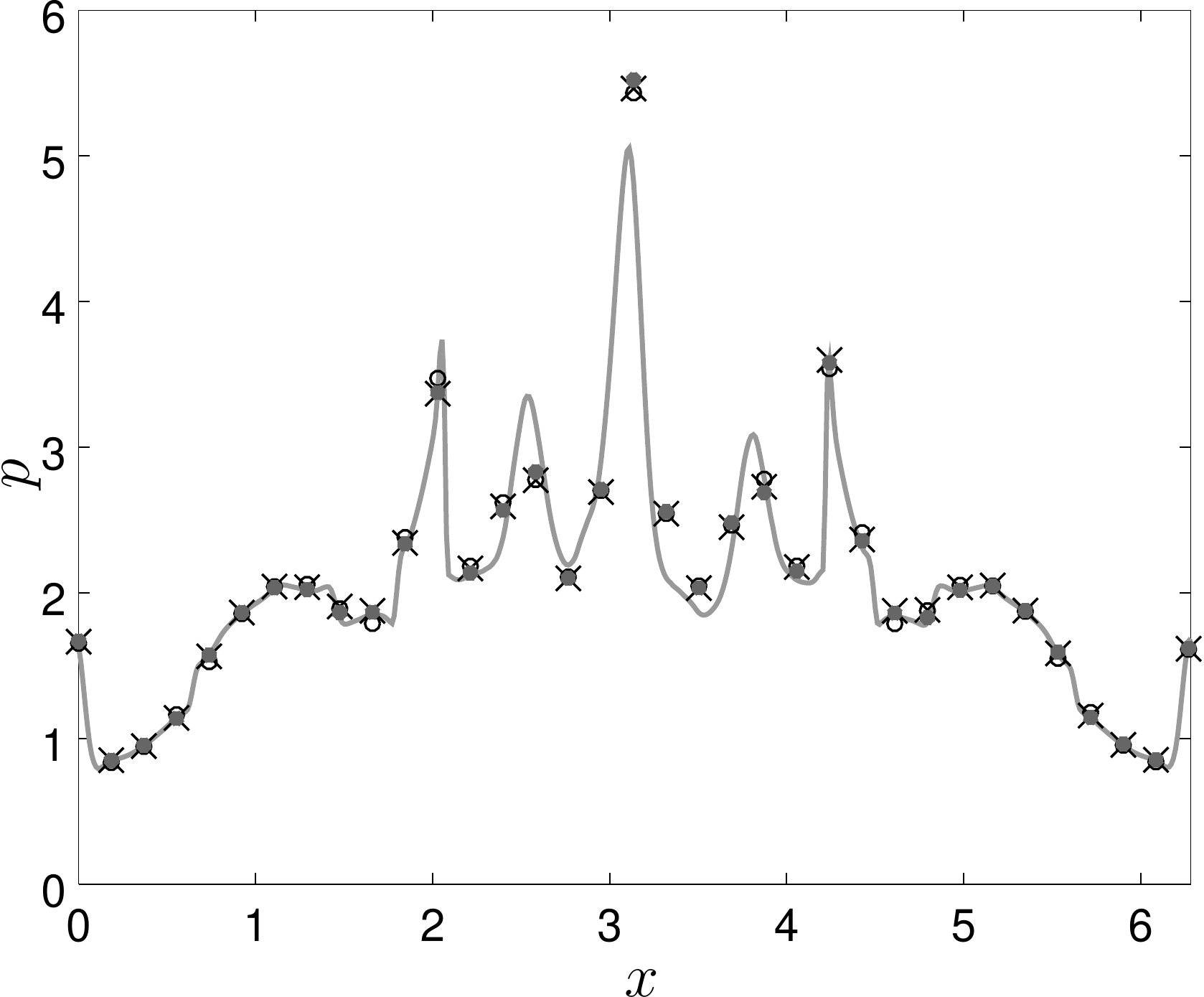}
\end{tabular}
\end{center}
\label{fig:cutsOTideal1}
\end{figure}

In the adaptive context, the simulations are performed by choosing fixed truncation parameters $\epsilon=0.01$, $\epsilon=0.03$ and the level-based one, $\epsilon_0=0.02$. In Figure~\ref{fig:OT2DL93}, we show the visualization of the MHD variables $p$ and $u_y$ obtained with $\epsilon=0.03$ at $t=\pi$. The symmetry of the solution is kept and the physical structures are well represented over the computational domain. The one dimensional cuts for $p$ at $x=\pi$ and $y=\pi$, including the reference solution, are presented in Figure~\ref{fig:cutsOTideal1}. We can notice that the similarity between the solutions, and observe the convergence towards the reference thus reproducing the expected physical behavior. 
In particular, the total pressure is a macroscopic entity that is the result of  the environment variables $\textbf{B}$, $\textbf{u}$ and $\rho$, which makes it a {\color{black} suitable} numerical sensor in studies such as instabilities, environment morphology, and reconnection phenomena.


\begin{figure}[H]
\psfrag{x}{$x$}
\psfrag{y}{$y$}
\caption{Adaptive meshes at $t=\pi$ and $L=9$, for the 2D Orszag--Tang vortex with $\epsilon = 0.01$ and $\epsilon=0.03$.}
\begin{center}
\begin{tabular}{cc}
$\epsilon=0.01$ & $\epsilon=0.03$\\
\includegraphics[width=0.4\linewidth]{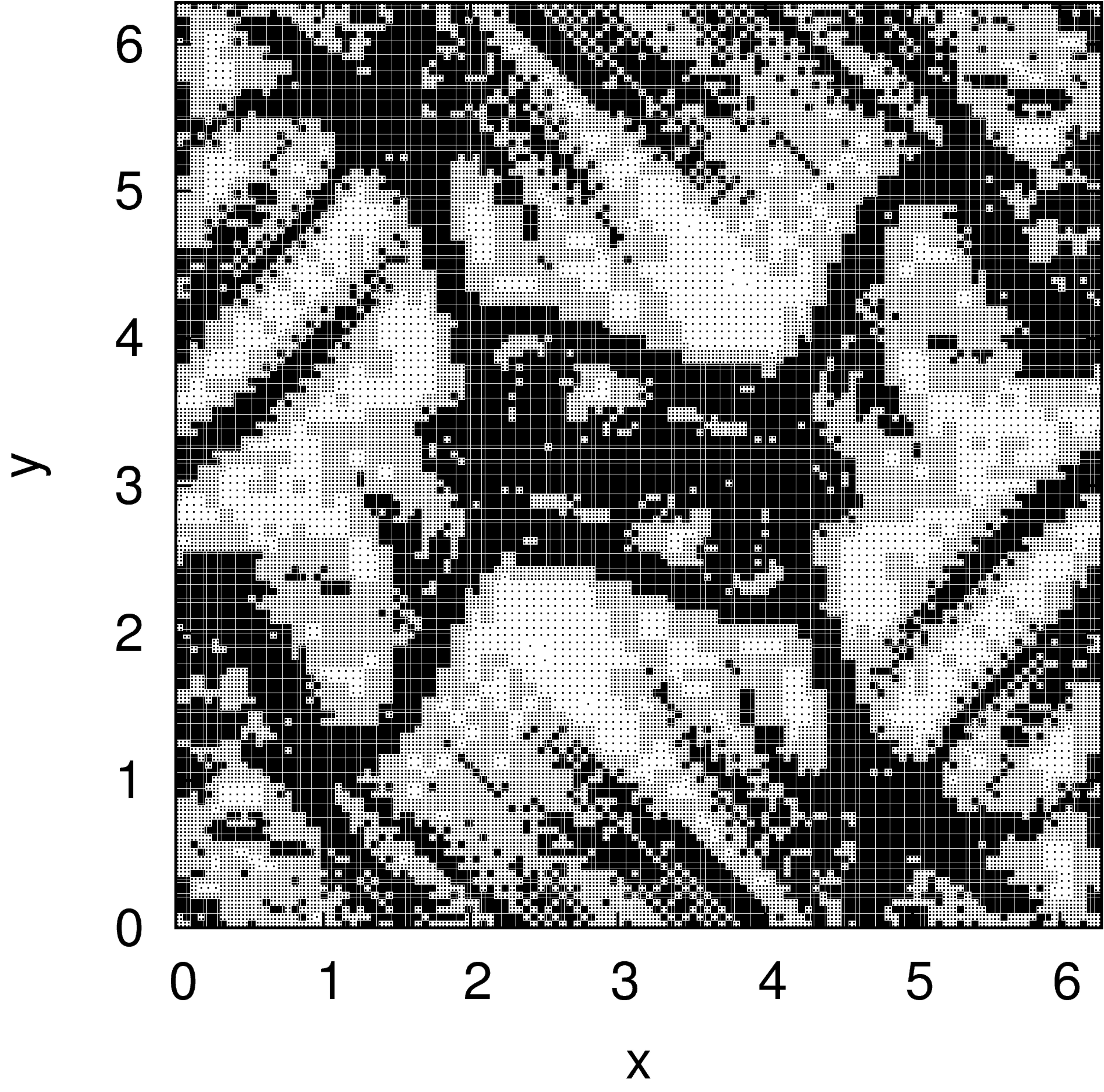} & \includegraphics[width=0.4\linewidth]{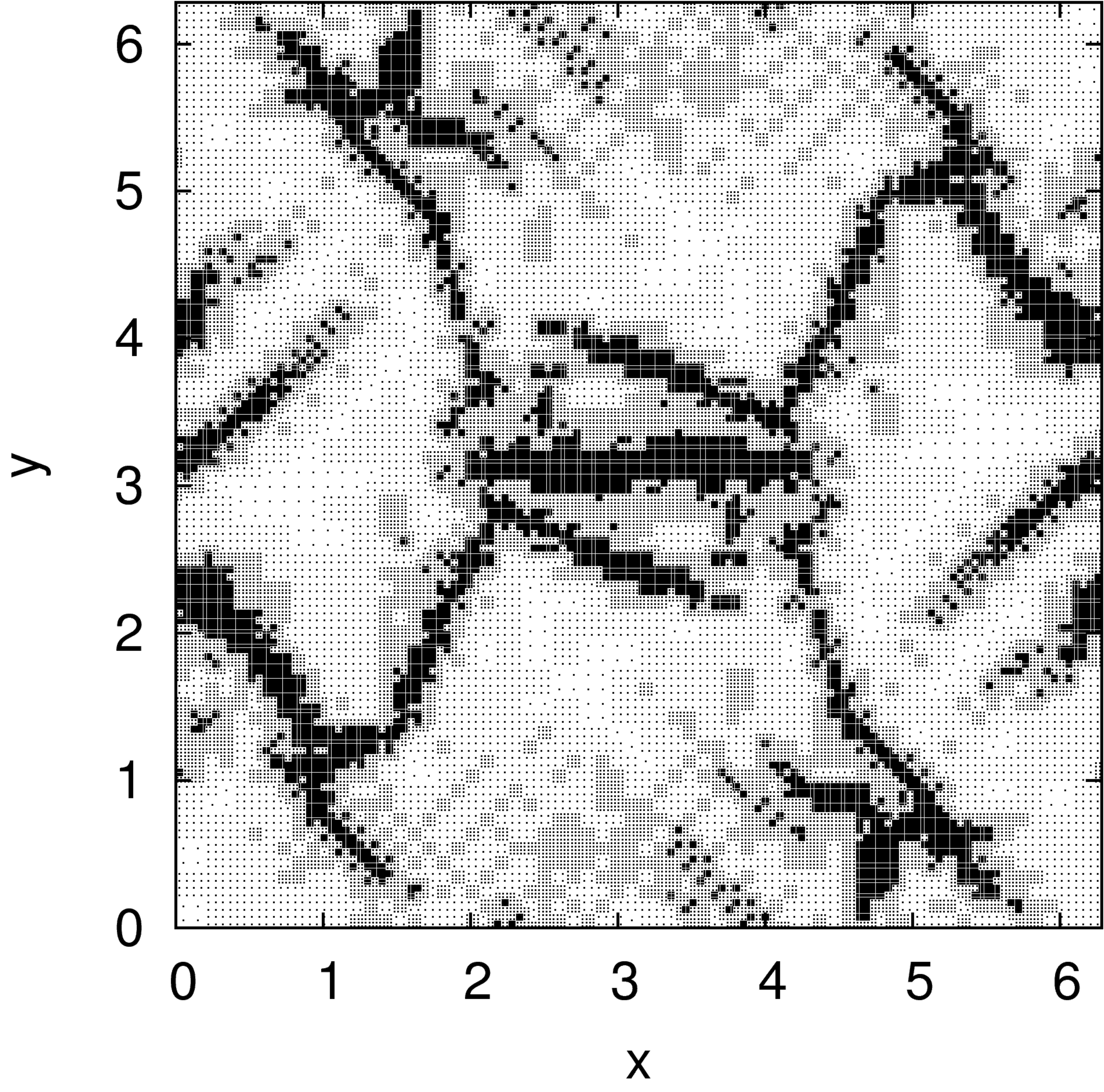}
\end{tabular}
\end{center}
\label{fig:OT2Dmesh}
\end{figure}

\begin{figure}[H]
\psfrag{L}{$L$}
\psfrag{EPS}{$\epsilon$}
\psfrag{L1}{$\mathcal{L}^1$}
\psfrag{L2}{$\mathcal{L}^2$}
\psfrag{MEM}{$\%$ célls}
\psfrag{CPU}{CPU Time}
\caption{Mass density $\rho$ errors at $t=\pi$ and $L=9$, for uniform and adaptive cases} 
\centering
\begin{tabular}{cc}
(a) $\epsilon$ $\times$ $\mathcal{L}^1$ error & (b) $\epsilon$ $\times$ $\mathcal{L}^2$ error \\
\includegraphics[width=0.4\linewidth]{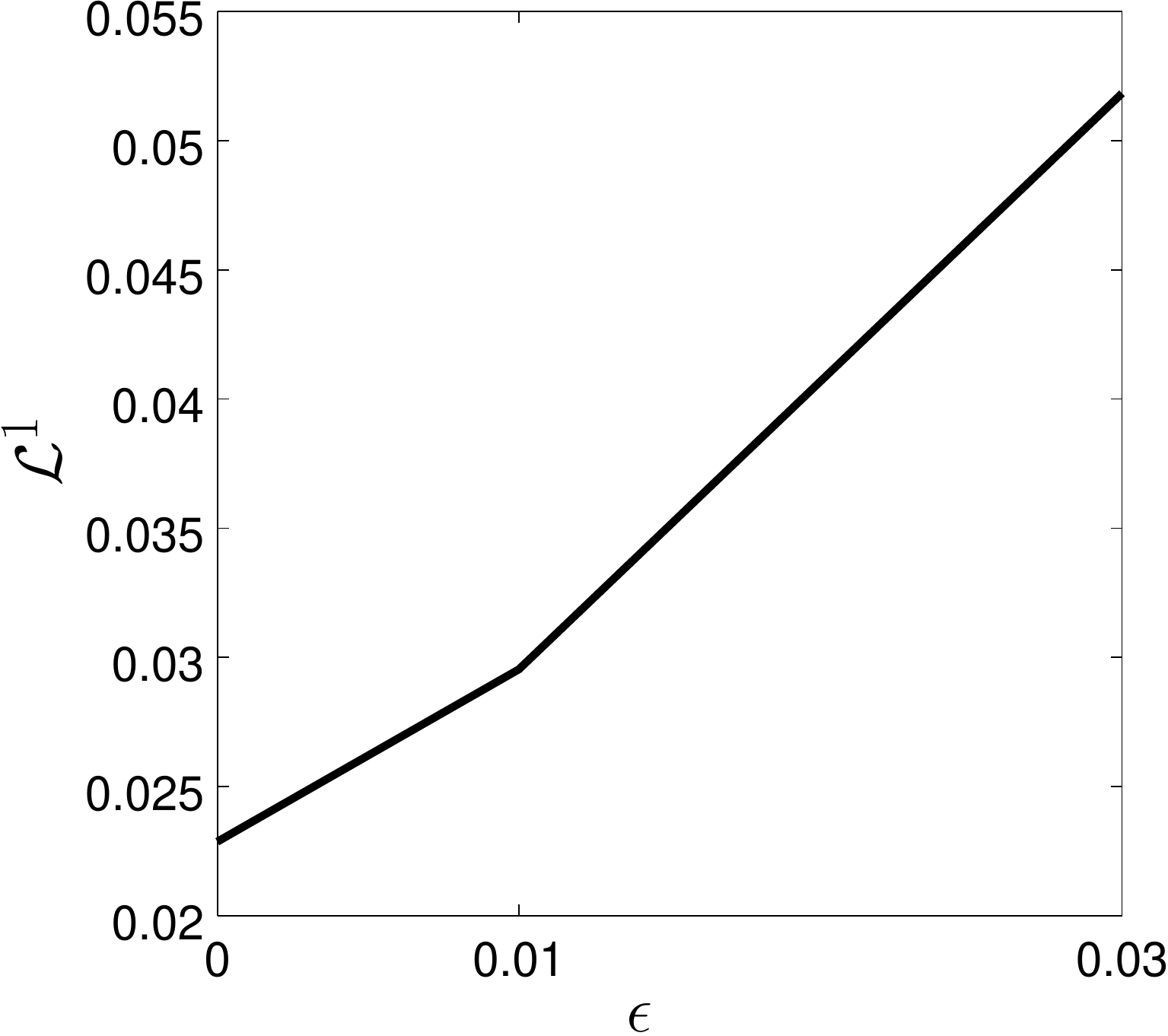} & \includegraphics[width=0.39\linewidth]{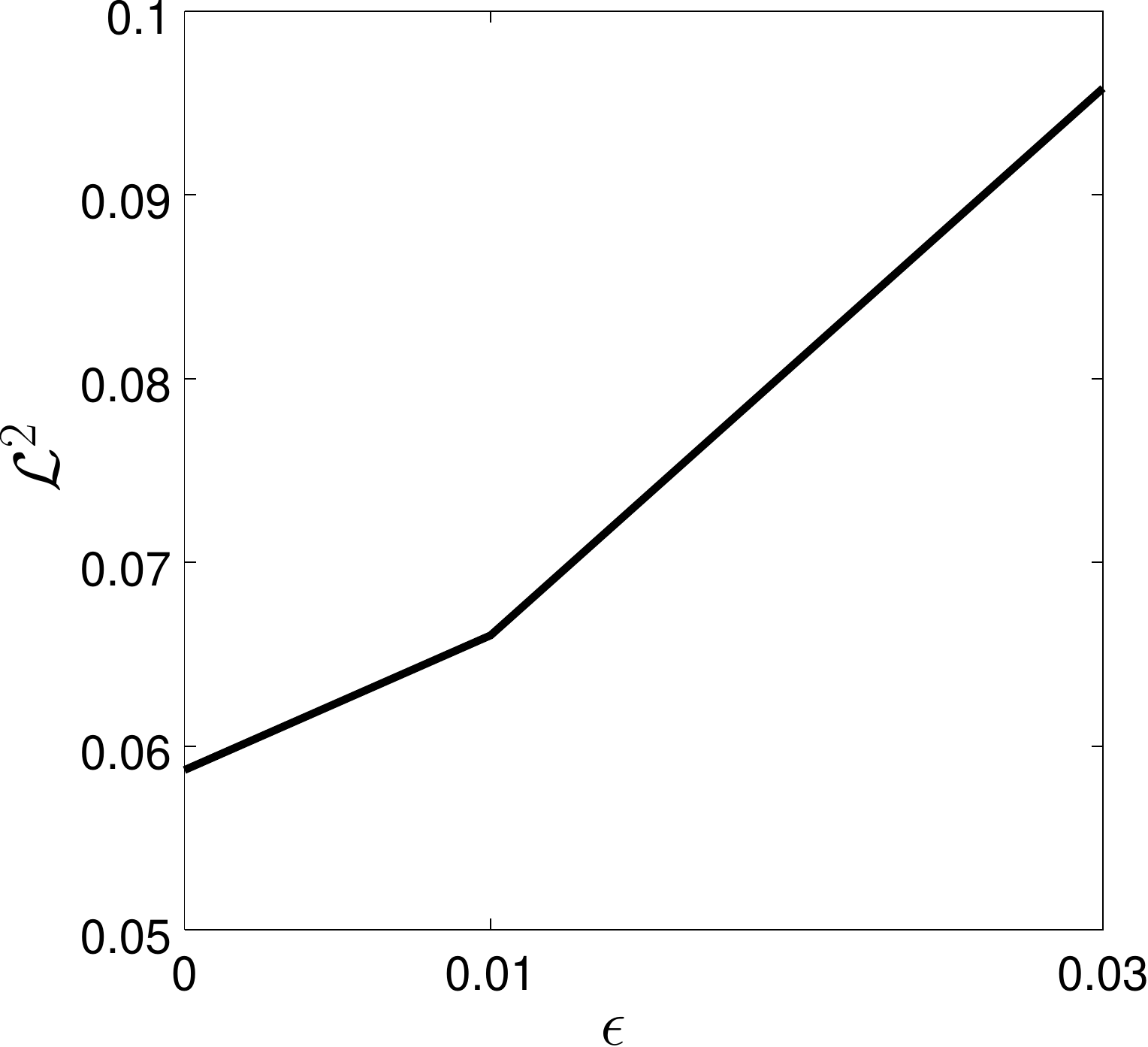}
\end{tabular}
\label{fig:epsxerr}
\end{figure}

The optimal threshold parameter found for this problem is $\epsilon=0.03$, because it provides a significant economy of cells and CPU time. Moreover, it also maintains the error in the same order of accuracy when compared to the other adaptive cases, as presented in Figure~\ref{fig:epsxerr} for the density variable. For $\epsilon=0.03$, only $26\%$ of the cells over time are needed for the simulation, decreasing the CPU time by $77\%$. The percentage of cells required for $\epsilon=0.01$ and $\epsilon_0=0.2$ simulations are, respectively, $45\%$ and $85\%$, leading to a $56\%$ and $11\%$ CPU time reduction. The adaptive meshes for $\epsilon=0.01$ and $\epsilon=0.03$ at $t=\pi$ are shown in Figure~\ref{fig:OT2Dmesh}. In both cases, the cells are located in regions that present stronger discontinuities, according to the visualization of the variables provided, allowing the structures to be well represented even in the case with less cells.

The $\mathcal{L}^1$ and $\mathcal{L}^2$ errors for the uniform and adaptive simulations are shown in Table~\ref{tab:ot2d-ideal-FV} for $p$ and $u_y$. 
The uniform mesh errors are slightly smaller, which is expected since the number of cells is significantly larger. In Figure~\ref{fig:epsxerr}, we show that as we increase the $\epsilon$ value, the error values also increase. We are comparing the uniform mesh, denoted by $\epsilon=0$, and two adaptive cases, $\epsilon=0.01$ and $\epsilon=0.03$. It is important to notice here that the $\epsilon$ value is also related to the number of cells, i.e., the number of cells tends to decrease as we increase $\epsilon$. Thus having more cells implies smaller. 
This shows that in the adaptive simulations we should find the optimal relation between computational gain and accuracy. We can also observe that the errors of the pressure variable are larger, which happens because this variable is obtained from the other $7$ MHD variables and thus their errors accumulate.

\begin{table}[H]
        \centering
        \caption{Errors for the ideal 2D O--T vortex for the uniform and adaptive, compared to the reference solution at $L=9$.}
        \vspace{2mm}
        \small{
        \begin{tabular}{@{}ccrr@{}}
        \toprule
        CARMEN--MHD &\multirow{2}{*}{Variables}       & \multicolumn{2}{c}{Errors $(\times 10^{-2})$} \\ \cmidrule{3-4}
        solver &&   $\mathcal{L}^1$     &    $\mathcal{L}^2$   
        \\ \cmidrule{1-4}
            \multirow{2}{*}{Uniform}
 		   	& $p   $ & 2.256  & 7.052 
 		   	\\
  			& $u_y $ & 0.628  & 1.653 
  			\\
  			 \multirow{1}{*}{Adaptive}            
 		  	& $p   $ & 5.337 &  11.79 
 		   	\\
  			($\epsilon=0.03$)  & $u_y $ & 1.954 &  3.582 
  			\\
  		\bottomrule
        \end{tabular}
        }
        \label{tab:ot2d-ideal-FV}
    \end{table}


\begin{figure}[H]
\psfrag{logEFV}{\small $log_{10}(\mathcal{E}_{MR}/\mathcal{E}_{UM})$}
\psfrag{t}{$t$}
\psfrag{VF}{UM}
\psfrag{eps001-46}{\small $\epsilon=0.01$}
\psfrag{eps02J-86}{\small $\epsilon_0=0.2$}
\psfrag{eps003-26}{\small $\epsilon=0.03$}
\caption{Global values of $\log_{10}(\mathcal{E}_{MR}/\mathcal{E}_{UM})$ over $t$, obtained for the uniform mesh (black) and adaptive cases with $\epsilon_0=0.2$ (circle), $\epsilon=0.01$ (cross) e $\epsilon=0.03$ (dotted).}
\begin{center}
\includegraphics[width=0.7\textwidth]{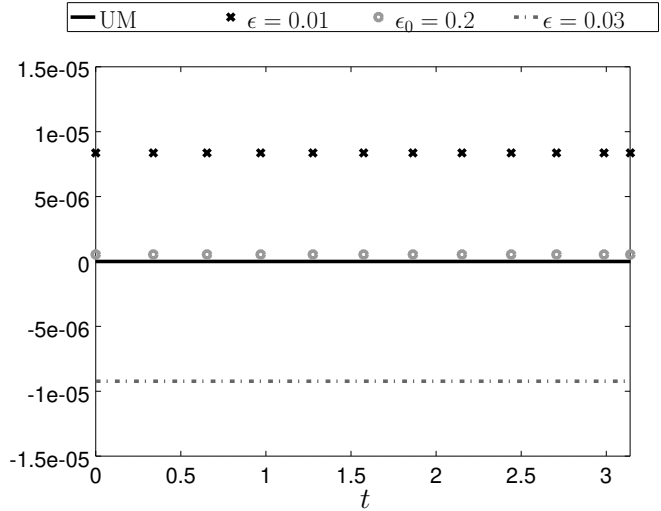}
\end{center}
\label{fig:enerPlot}
\end{figure}

\begin{figure}[htb!]
\psfrag{EPS}{$\epsilon$}
\psfrag{DIVERROR}{\small{$\max\left(d V \frac{|\nabla\cdot\textbf{B}|}{|\textbf{B}|}\right)$}}
\caption{Quantity $\epsilon$ $\times$ $\varepsilon_{div}$ at $t=\pi$ and $L=9$, for uniform ($\epsilon=0$) and adaptive cases with $\epsilon=0.01$, $\epsilon=0.03$.}
\centering
\includegraphics[width=0.5\linewidth]{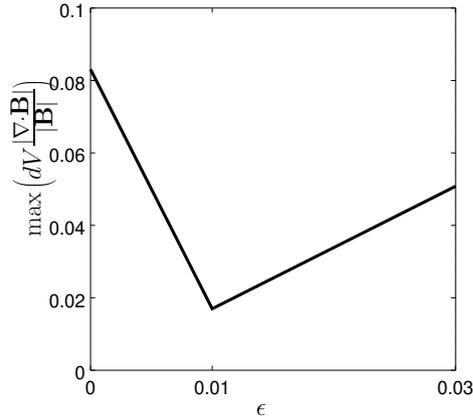}
\label{fig:epsxdiv}
\end{figure}

To evaluate the conservation of the total energy density $\mathcal{E}$, we present the values referring to $\log_{10}(\mathcal{E}_{MR}/\mathcal{E}_{UM})$ in Figure~\ref{fig:enerPlot}, where $\mathcal{E}_{MR}$ and $\mathcal{E}_{UM}$ are the integral values of $\mathcal{E}$ over time on adaptive and uniform meshes, respectively. By using this measure, it is possible to study the energy conservation for $\epsilon_0=0.2$ (circle), $\epsilon=0.01$ (cross) and $\epsilon=0.03$ (dotted) and also to verify how close these values are to the uniform case (solid line). The conservation of energy over time holds for every case presented, maintaining the physical properties. As much as we increase the number of cells of the simulations, the integral values converge to the uniform case, \textit{e.g.} the $\epsilon_0=0.2$ case.
In Figure~\ref{fig:epsxdiv}, values of the divergence error $\varepsilon_{div}$ are presented for two different values of $\epsilon$ and the uniform mesh computation with $L=9$. For all cases, the values are below $10^{-1}$ and, consequently, the computations satisfy the restriction $\varepsilon_{div}<1$. 
This shows that the GLM-MHD parabolic-hyperbolic strategy combined with the adaptive MR technique still presents small numerical values of $\nabla\cdot\textbf{B}$, as expected. Hence this combination preserves the desired precision of the numerical solution at the final time of the simulation.

\subsection{Shock-cloud interaction}
\label{subsec:cloud}
The shock-cloud interaction models the disruption of a dense cloud with a shock-wave. The problem presents strong discontinuities and it is a challenging test case to evaluate the stability of numerical schemes. This type of simulation was firstly proposed in \cite{dai1998simple} and the initial condition we use is based on \cite{Touma2006617}, given in Table~\ref{tab:CISC} and with $u_y=u_z=B_x=0.0$. We consider a circular cloud, 10 times denser than its background, with radius $r=0.15$ and center $(0.25,0.5,0.5)$ in three dimensions and $(0.25,0.5)$ in two dimensions. The shock is located at $x=0.05$.

As simulation parameters we choose the physical time $t=0.06$, the divergence cleaning parameter $\alpha= 0.4$, the threshold parameter $\epsilon^0 =0.01$ and $\gamma = 5/3$.

\begin{figure}[H]
\caption{Variables $\rho$ and $B_z$ at $t=0.06$ and $L=9$, obtained with the CARMEN--MHD code for ideal MHD, for the 2D shock-cloud problem and $\epsilon_0=0.01$.}
\begin{center}
\begin{tabular}{cc}
$\rho$ & $B_z$\\
\includegraphics[width=0.44\textwidth]{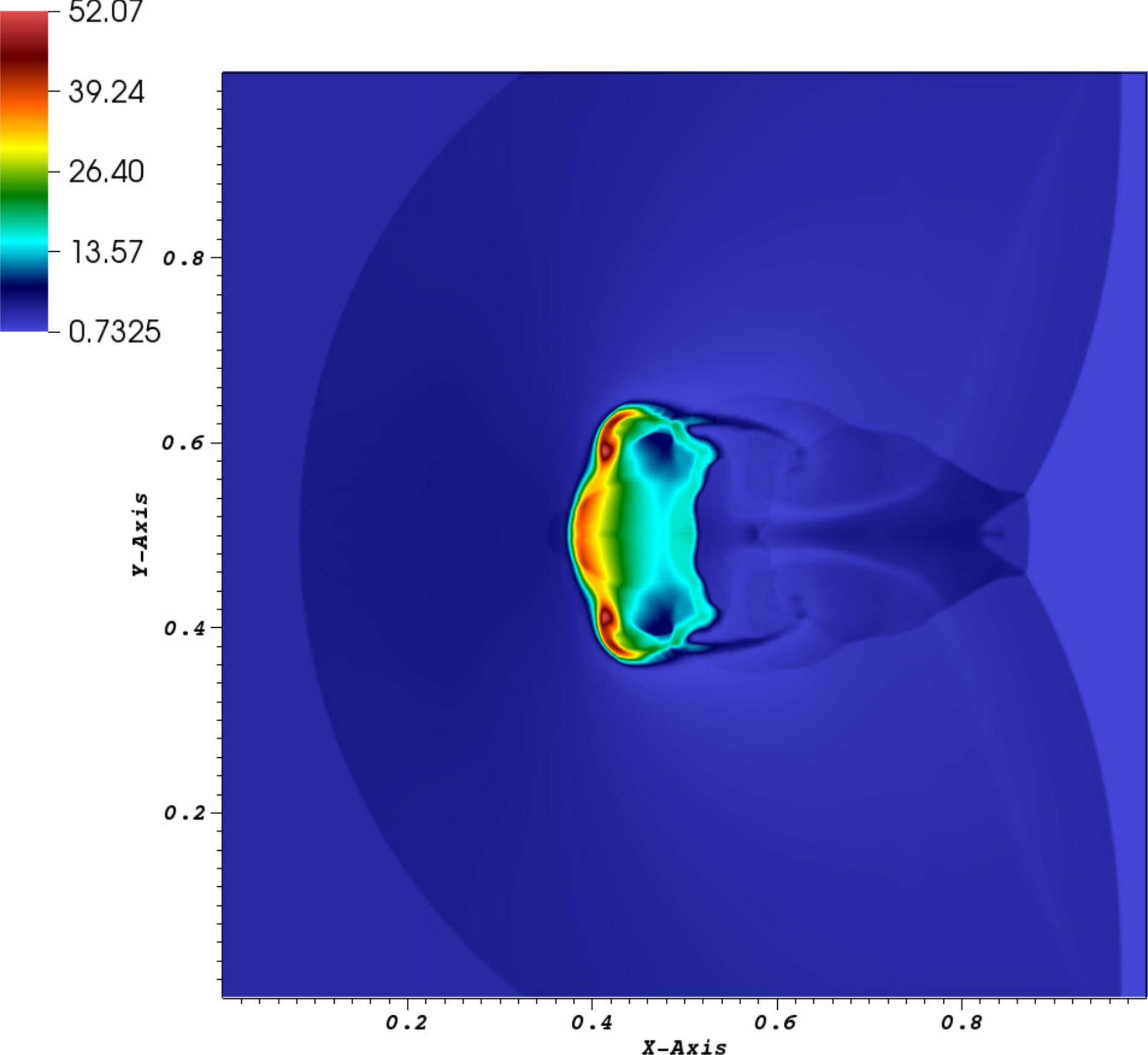}&\includegraphics[width=0.44\textwidth]{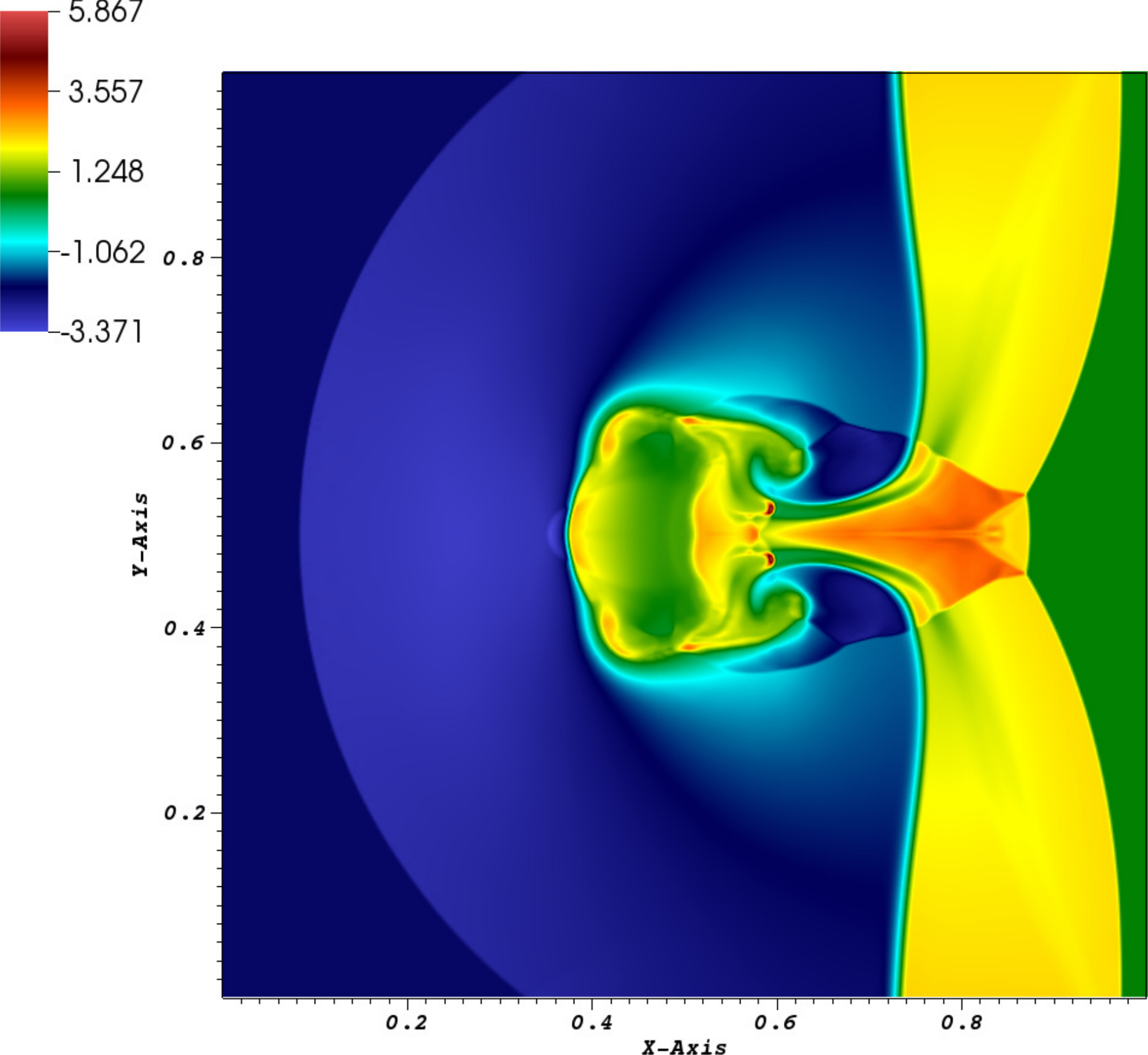}  
\end{tabular}
\end{center}
\label{fig:SC2D1}
\end{figure}

\begin{table}[H]
    \centering{
    \caption{Shock-cloud initial condition.}
        \begin{tabular}{@{}cccccc@{}}
        \toprule
         &$\rho$ & $p$ & $u_x$  & $B_y$ & $B_z$ \\
        \cmidrule{2-6} 
             $x\leq 0.05$ & 3.86859000 & 167.34500000 & \phantom{1}0.00000000 & 2.18261820 & -2.18261820 \\
             $x>0.05$     & 1.00000000 & \phantom{67}1.00000000 & 11.25360000 & 0.56418958 & \phantom{-}0.56418958 \\[2mm]
        \bottomrule 
        \end{tabular}}
    \label{tab:CISC}
\end{table}
In two dimensions, we consider a high density cloud such as a circle centered in $(0.25,0.5)$, with radius $r=0.15$ and $\rho=10$.
The variables $\rho$ and $B_z$ are presented in Figure~\ref{fig:SC2D1}, obtained with the CARMEN--MHD simulations. We can notice there is a strong discontinuity present in the interval $x\in[0.4,0.5]$, which appears after the explosion coming from the interaction between the shock and the cloud. There are several local structures that are well represented in the visualization.

\begin{figure}[H]
\psfrag{RHO}{$\rho$}
\psfrag{VX}{$u_x$}
\psfrag{BZ}{$B_z$}
\psfrag{y}{$y$}
\psfrag{x}{$x$}
\psfrag{EPS001}{\scalebox{.5}{$\epsilon_0=0.01$}}
\caption{Cuts of the variable $B_z$  at $t=0.06$,  and $L=9$, for reference solution (solid line) and CARMEN--MHD adaptive solutions with $\epsilon_0 = 0.01$ (cross).}
\begin{center}
\begin{tabular}{cc}
$x=0.5$ & $y=0.5$ \\
\includegraphics[width=0.45\linewidth]{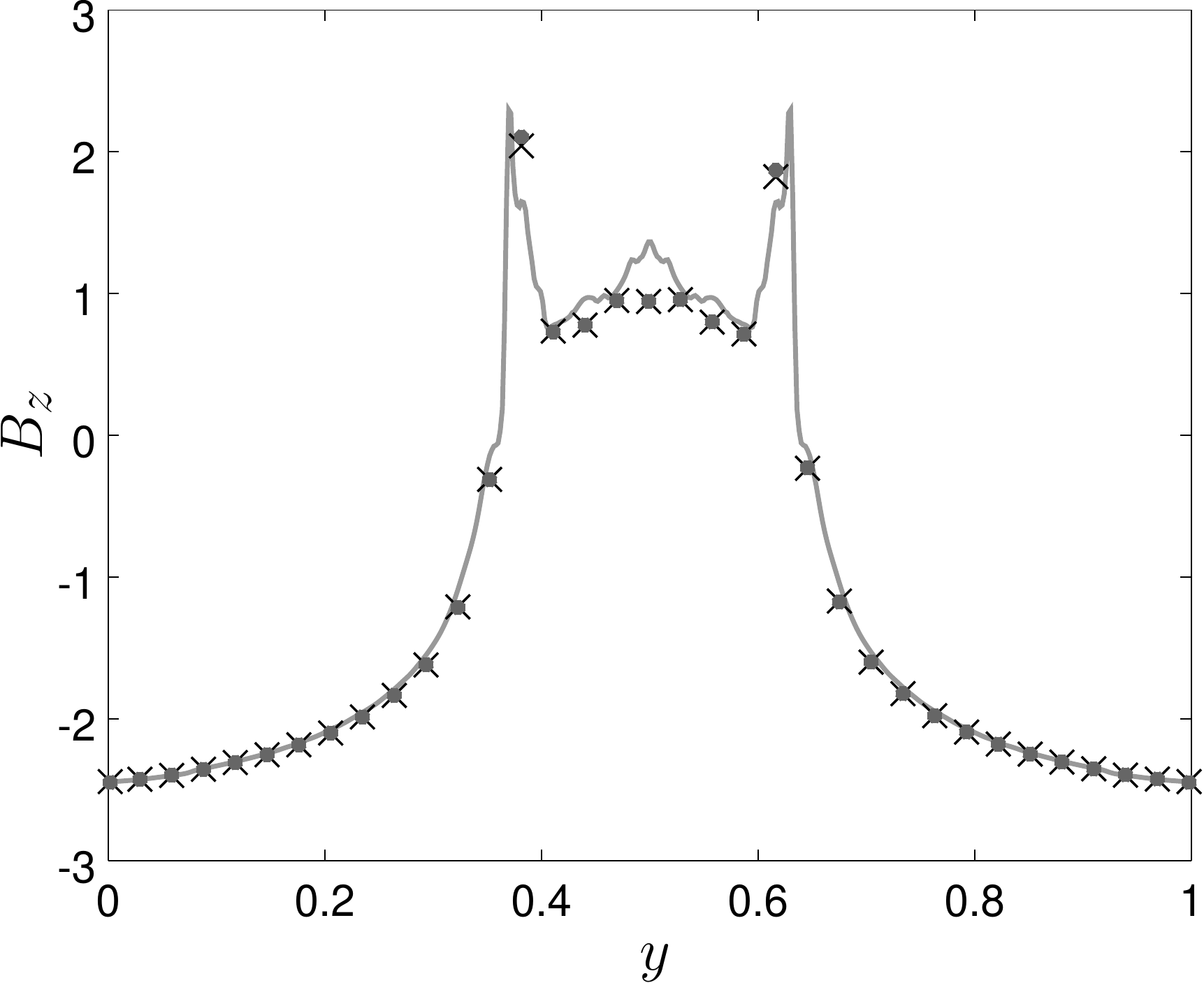} & \includegraphics[width=0.45\linewidth]{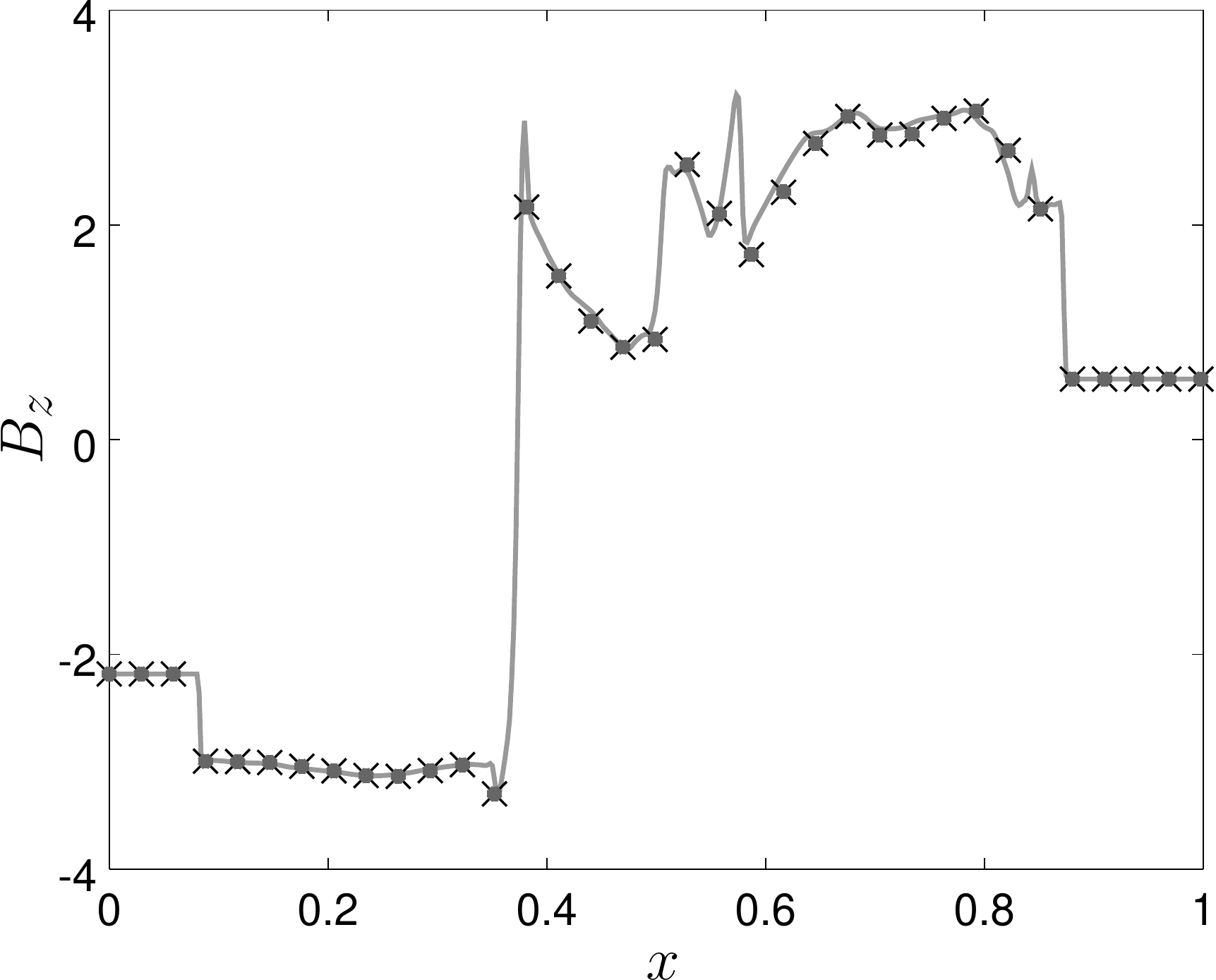}
\end{tabular}
\end{center}
\label{fig:cutsSCideal}
\end{figure}

In Figure~\ref{fig:cutsSCideal}, we present cuts of the variable $B_z$ at $x=0.5$, $y=0.5$ and $t=0.06$, for $\epsilon_0=0.01$ and $\epsilon=0.003$. The variable $B_z$ presents many discontinuities over the domain and the adaptive MR approach captures them. Moreover, these results reinforce the convergence of our numerical solutions. A slightly different topology may be observed in the solution, mainly due to the different numerical schemes. Nevertheless the solutions present the same behavior and {\color{black} have the same} accuracy order.

\begin{figure}[H]
\caption{Adaptive meshes at $t=0$, $t=0.06$ and $L=9$, for the ideal 2D shock-cloud problem with $\epsilon_0 = 0.01$.}
\begin{center}
\begin{tabular}{cc}
$t=0$ & $t=0.06$\\
\includegraphics[width=0.4\linewidth]{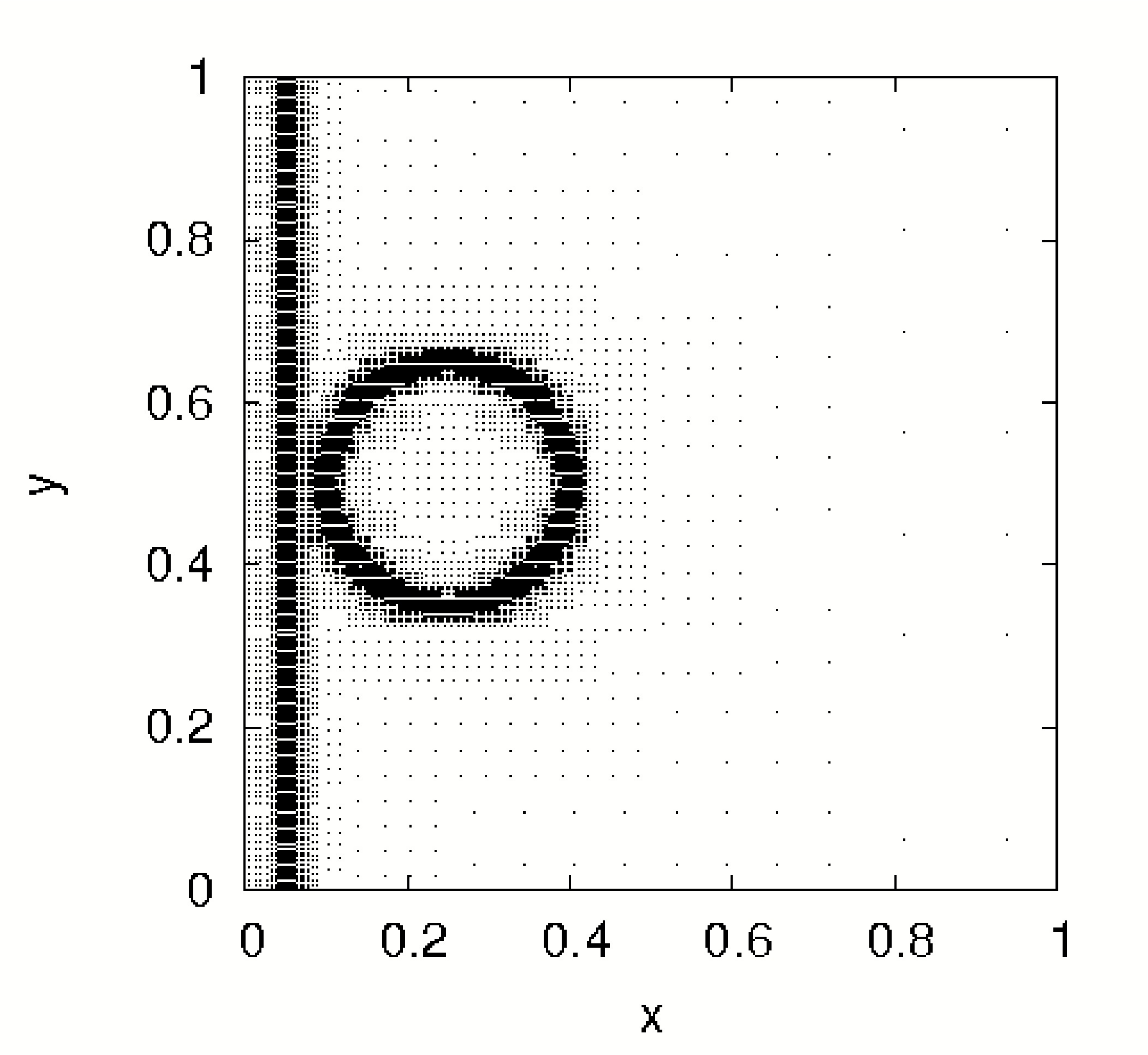} & \includegraphics[width=0.4\linewidth]{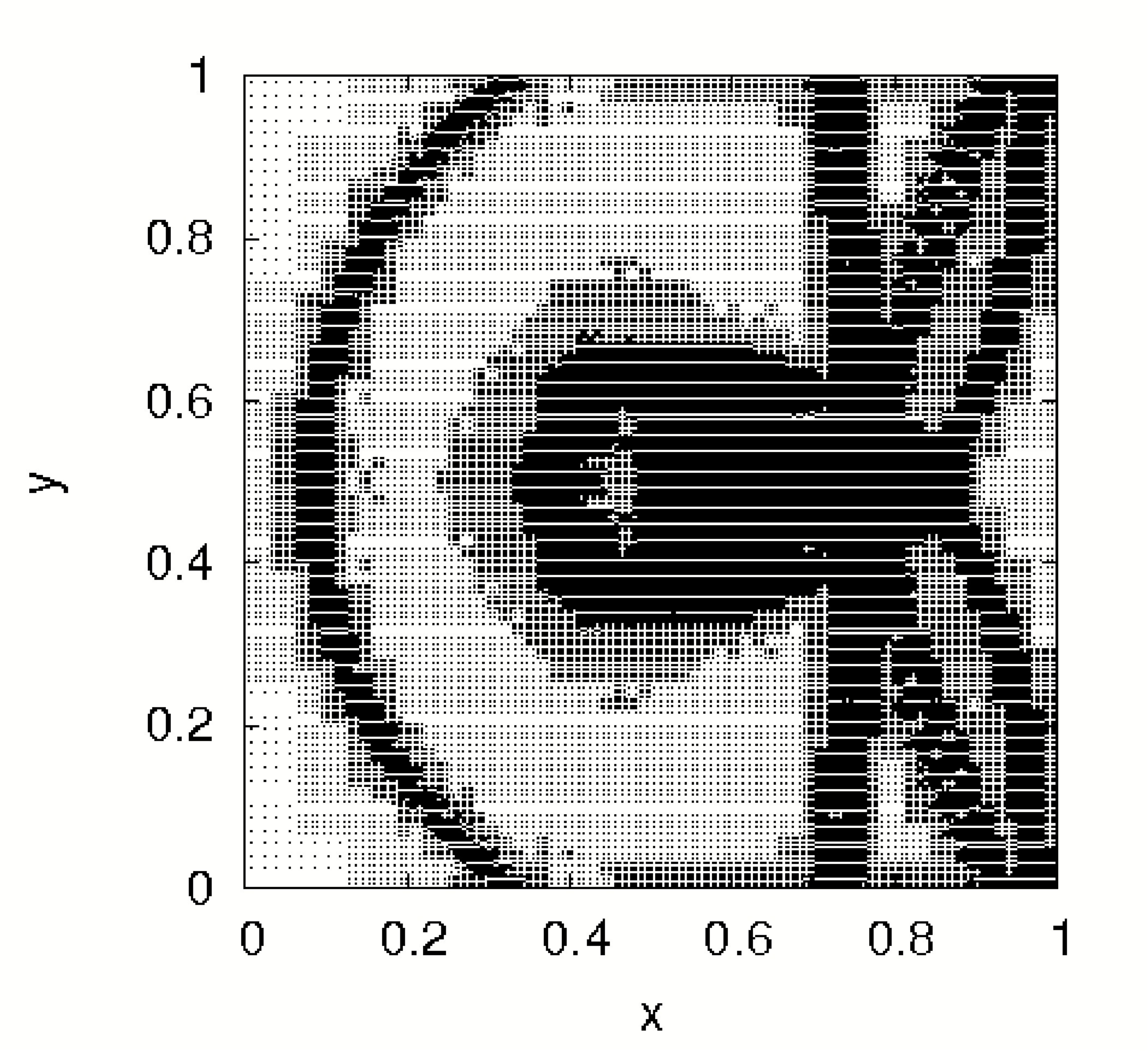}
\end{tabular}
\end{center}
\label{fig:SC2Dmesh}
\end{figure}

The adaptive meshes are presented in Figure~\ref{fig:SC2Dmesh} for simulations with threshold parameter $\epsilon_0=0.01$, that uses approximately $40\%$ of the cells over time at $t=0.06$. It implies a $55\%$ reduction of CPU time. The meshes outline the local structures of the problem at times $t=0$ and $t=0.06$, with most of the cells being concentrated in the shock front and the border of the cloud for the initial time, and on the sharper local structures for the final time.
In Table~\ref{tab:sc2d-ideal}, the errors evaluate quantitatively the proximity of the results obtained with uniform and adaptive meshes.  The $\mathcal{L}^1$ and $\mathcal{L}^2$ errors remain of order $10^{-1}$ for both cases, showing that it is possible to obtain the same order of accuracy by using only $40\%$ of the cells. We should recall at this point that the complexity of the shock-cloud problem, as it models an explosion with strong discontinuities, can influence on the elevated error values presented.

\begin{table}[H]
        \centering
        \caption{Errors obtained for the ideal 2D shock-cloud problem compared to the reference solution at level $L=9$.}
        \vspace{2mm}
        \small{
        \begin{tabular}{@{}ccrrr@{}}
        \toprule
        CARMEN--MHD &\multirow{2}{*}{Variables}       & \multicolumn{2}{c}{Errors $(\times 10^{-1})$} \\\cmidrule{3-4}
                  solver &&   $\mathcal{L}^1$     &    $\mathcal{L}^2 $ 
                  \\ \cmidrule{1-4}
            \multirow{2}{*}{Uniform}
            	& $\rho$  & 1.543 &  11.55 \\ 
           	& $B_z $  & 0.391 &  1.509 \\
        \midrule
            \multirow{1}{*}{Adaptive}
            	& $\rho$ & 1.539  &  11.60 \\
   		($\epsilon^0=0.01$)	& $B_z $ & 0.390  &  1.511 \\ 
        \bottomrule
        \end{tabular}
        }
        \label{tab:sc2d-ideal}
    \end{table}

\begin{figure}[H]
\caption{Variables $\rho$ and $B_z$ at $t=0.06$ and $L=7$, obtained with ideal CARMEN--MHD code for the 3D shock-cloud problem and $\epsilon_0=0.01$.}
\begin{center}
\begin{tabular}{cc}
$\rho$ & $B_z$\\ 
  \includegraphics[width=0.44\textwidth]{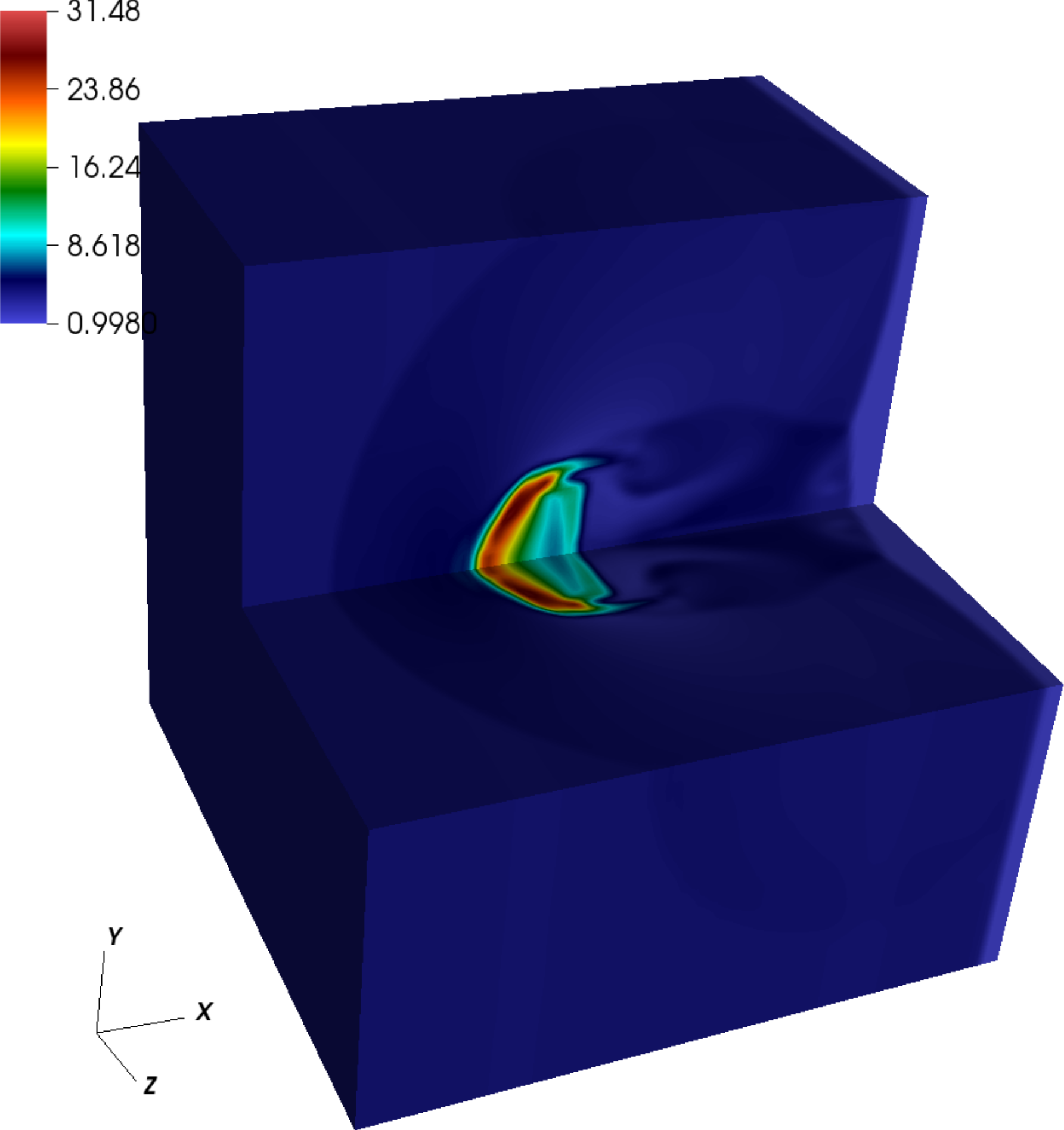}&\includegraphics[width=0.44\textwidth]{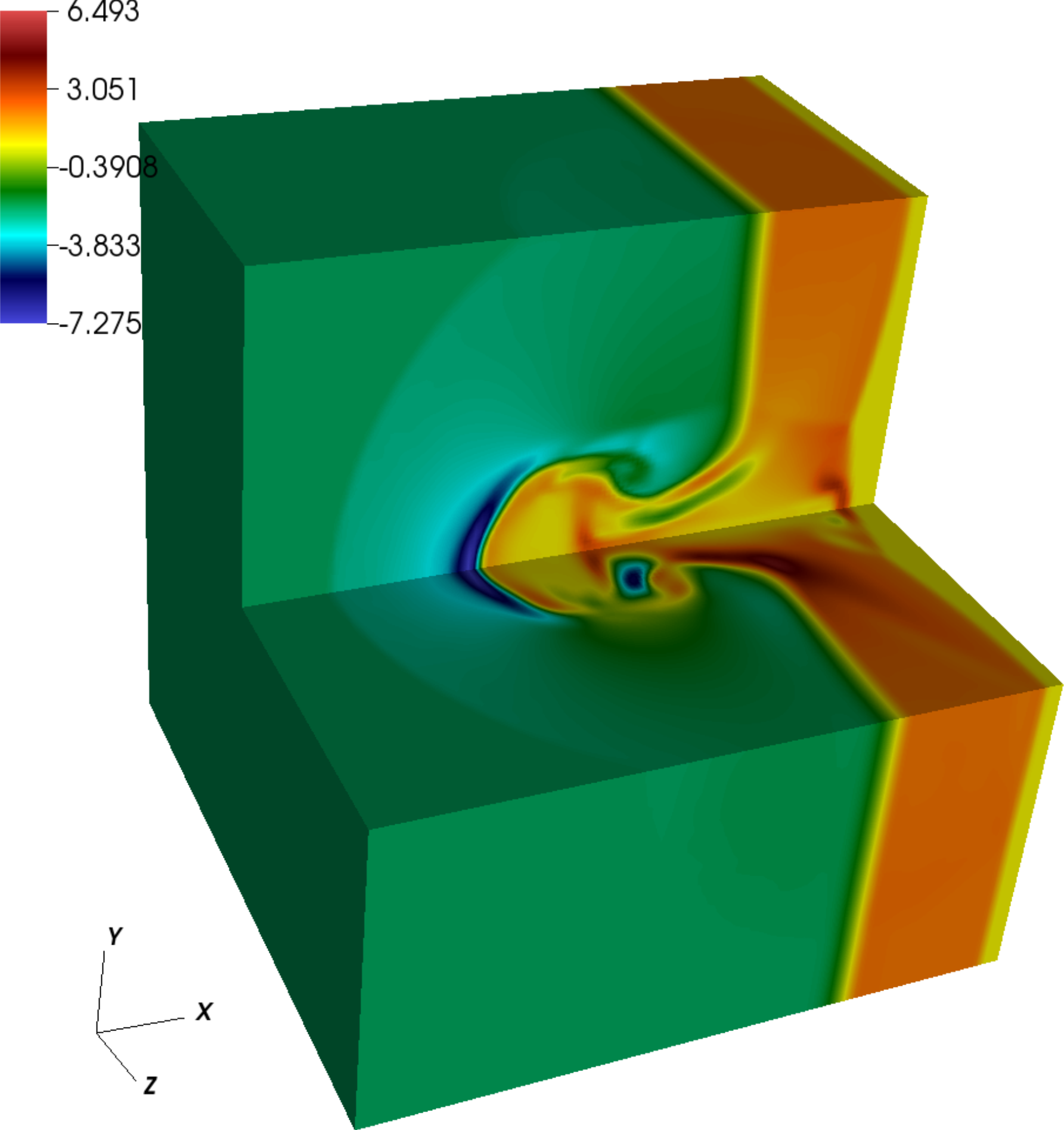} 
\end{tabular}
\end{center}
\label{fig:SC3D1}
\end{figure}

In three dimensions, the magnetic cloud is centered on $(0.25,0.5,0.5)$. We present the 3D shock-cloud simulation at level $L=7$ ($128^3$ cells) and $\epsilon_0=0.01$. Figure~\ref{fig:SC3D1} shows variables $\rho$ and $B_z$ obtained with the ideal MHD model at $t=0.06$. The adaptive structures are coherent with the expected ones from the reference. The approximation errors are presented in Table~\ref{tab:sc3d-ideal}. The errors $\mathcal{L}^1$ and $\mathcal{L}^2$ stay at the order $10^{-3}$ and $10^{-4}$, respectively, for both adaptive and uniform cases. If we compare these values to the two dimensional case, we can notice that they decrease for the three dimensional case. This probably occurs because of the $z$ component, which can present a solution with a globally smoother behavior compared to the two dimensional case, where complex structures are located all over the domain.

\begin{table}[H]
        \centering
        \caption{Errors obtained for the ideal 3D shock-cloud problem for the uniform and adaptive ($\epsilon^0=0.01$), compared to the reference solution at level $L=7$.}
        \vspace{2mm}
        \small{
        \begin{tabular}{@{}clrr@{}}
        \toprule
        CARMEN--MHD &\multirow{2}{*}{Variables}       & \multicolumn{2}{c}{Errors $(\times 10^{-3})$} \\\cmidrule{3-4}
           solver & &   $\mathcal{L}^1$     &    $\mathcal{L}^2$     \\ 
           \cmidrule{1-4}
            \multirow{2}{*}{Uniform}
			& $\rho$ & 0.143 &  4.798 \\ 
   			& $B_z $ & 0.132 &  3.966 \\ 
        \toprule
            \multirow{2}{*}{Adaptive}
           	& $\rho$ & 0.146 &  4.818 \\ 
   			& $B_z $ & 0.138 &  4.089  \\
        \bottomrule
        \end{tabular}
        }      
        \label{tab:sc3d-ideal}
\end{table}

\begin{figure}[H]
\caption{Ideal 3D shock-cloud problem: adaptive mesh ($\epsilon_0=0.01$) sections on planes $xy$, $yz$, $zx$ at $t=0.06$ and $L=7$.}
\begin{center}
\begin{tabular}{ccc}
$xy$ & $yz$ & $zx$\\
\includegraphics[width=0.3\textwidth]{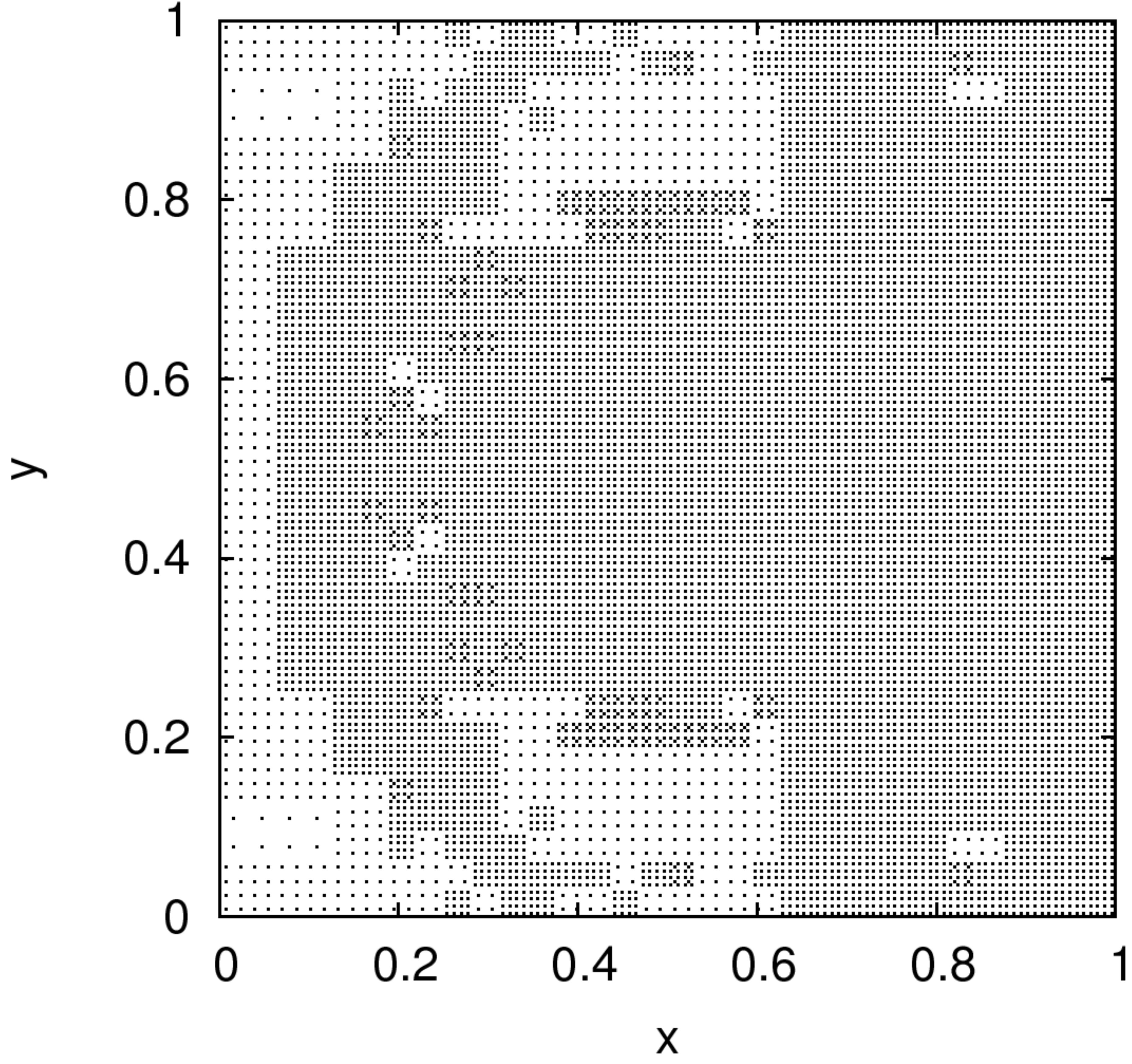} & \includegraphics[width=0.3\textwidth]{SC_3D_L7_yz-eps-converted-to.pdf} & \includegraphics[width=0.3\textwidth]{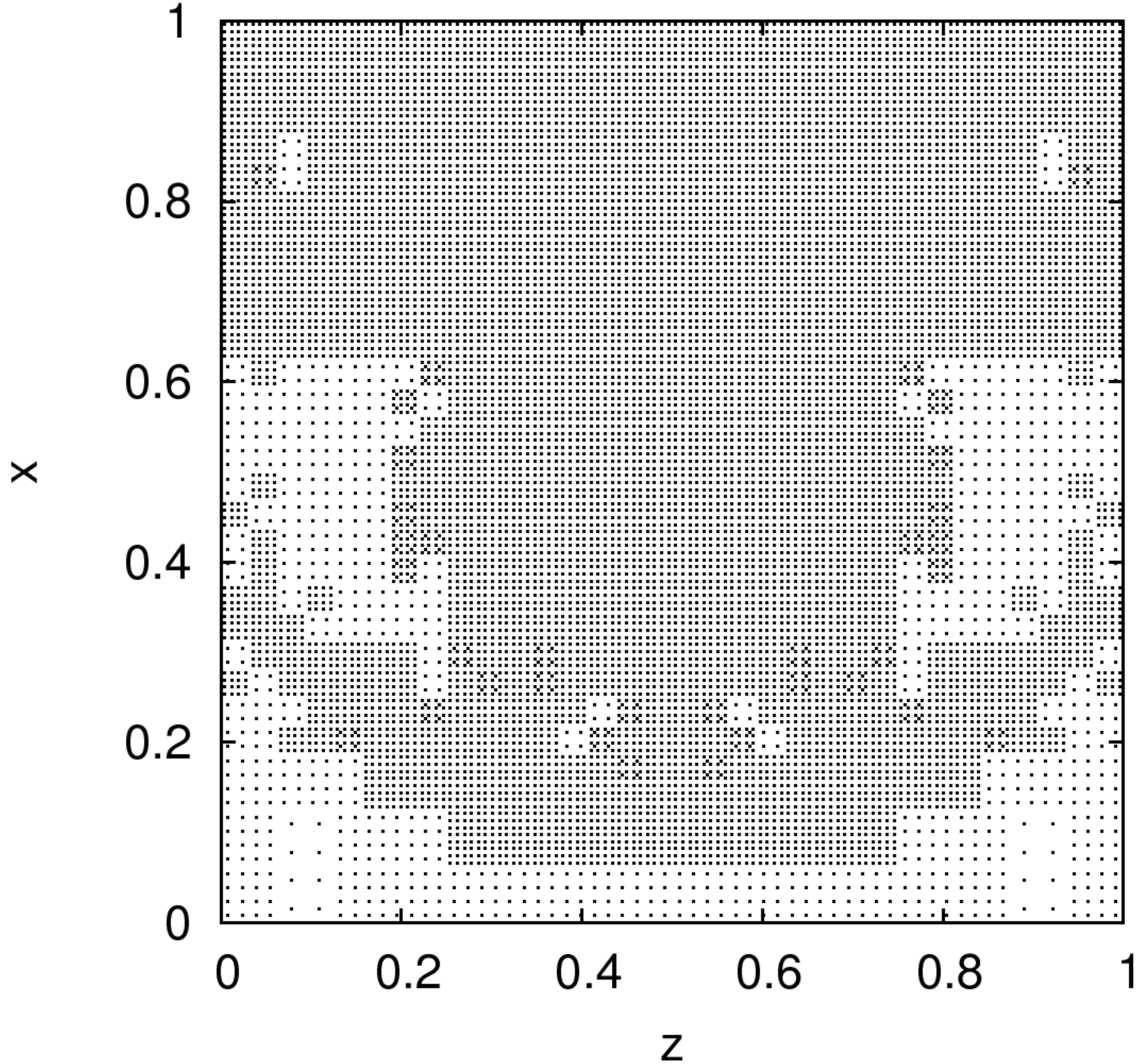}
\end{tabular}
\end{center}
\label{fig:SC3DL7mesh}
\end{figure}

The adaptive mesh is presented in Figure~\ref{fig:SC3DL7mesh} for $\epsilon_0=0.01$, decomposed according to the $xy$, $yz$ and $zx$, as we take sections at the interval $[0.45,0.55]$ on $z$, $x$ and $y$ axes, respectively. This type of visualization makes the mesh adaptivity clearer for the 3D case. It is possible to observe that the cells of this simulation are located exactly in regions of the stronger discontinuities. The $xy$ mesh is similar to the 2D case, in a coarser level. The 3D simulation demands $58\%$ of the cells over time, when compared to a uniform mesh, and causes a reduction of $34\%$ in CPU time. When we increase the refinement level to $L=8$, the percentage of cells {\color{black} required} decreases to $42\%$ and reduces the CPU time by $54\%$. This suggests that the adaptivity tends to improve as we increase the maximum level of refinement of the 3D case.

For the simulation of the resistive 3D shock-cloud, we add a constant resistivity $\eta=0.02$ all over the computational domain. 
Numerical tests showed that this value is reasonable for this problem, as the diffusive effect is sufficient to smooth the local structures of the problem without losing its intrinsic topology. 

\begin{figure}[H]
\caption{
Constant resistive 3D shock-cloud problem: adaptive simulation for variables $\rho$ and $B_z$ with $\epsilon_0=0.01$  at $t=0.06$ with $L=6$, and $\eta=0.02$}
\begin{center}
\begin{tabular}{cc}
$\rho$ & $B_z$\\
\includegraphics[width=0.44\textwidth]{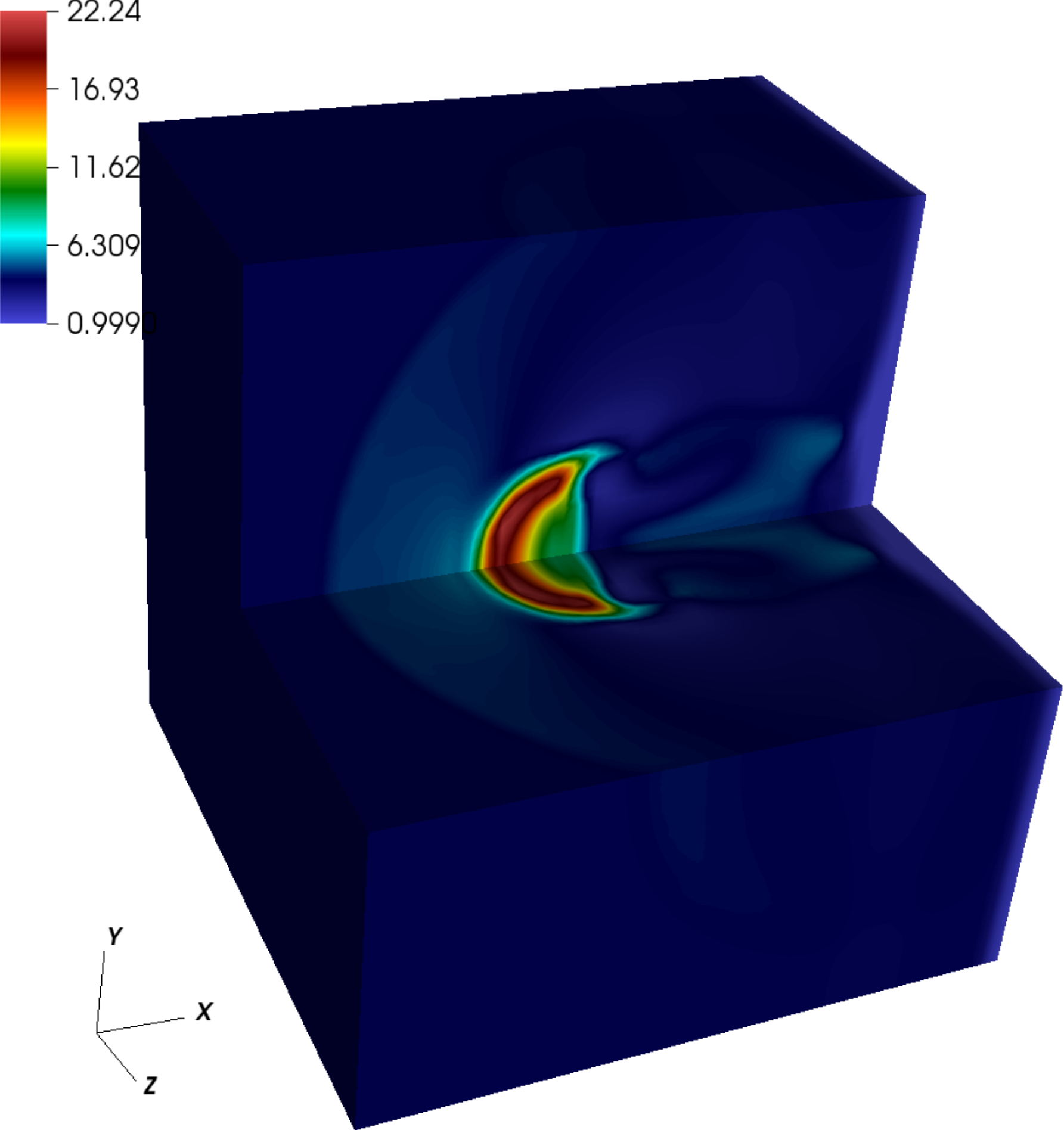}&\includegraphics[width=0.44\textwidth]{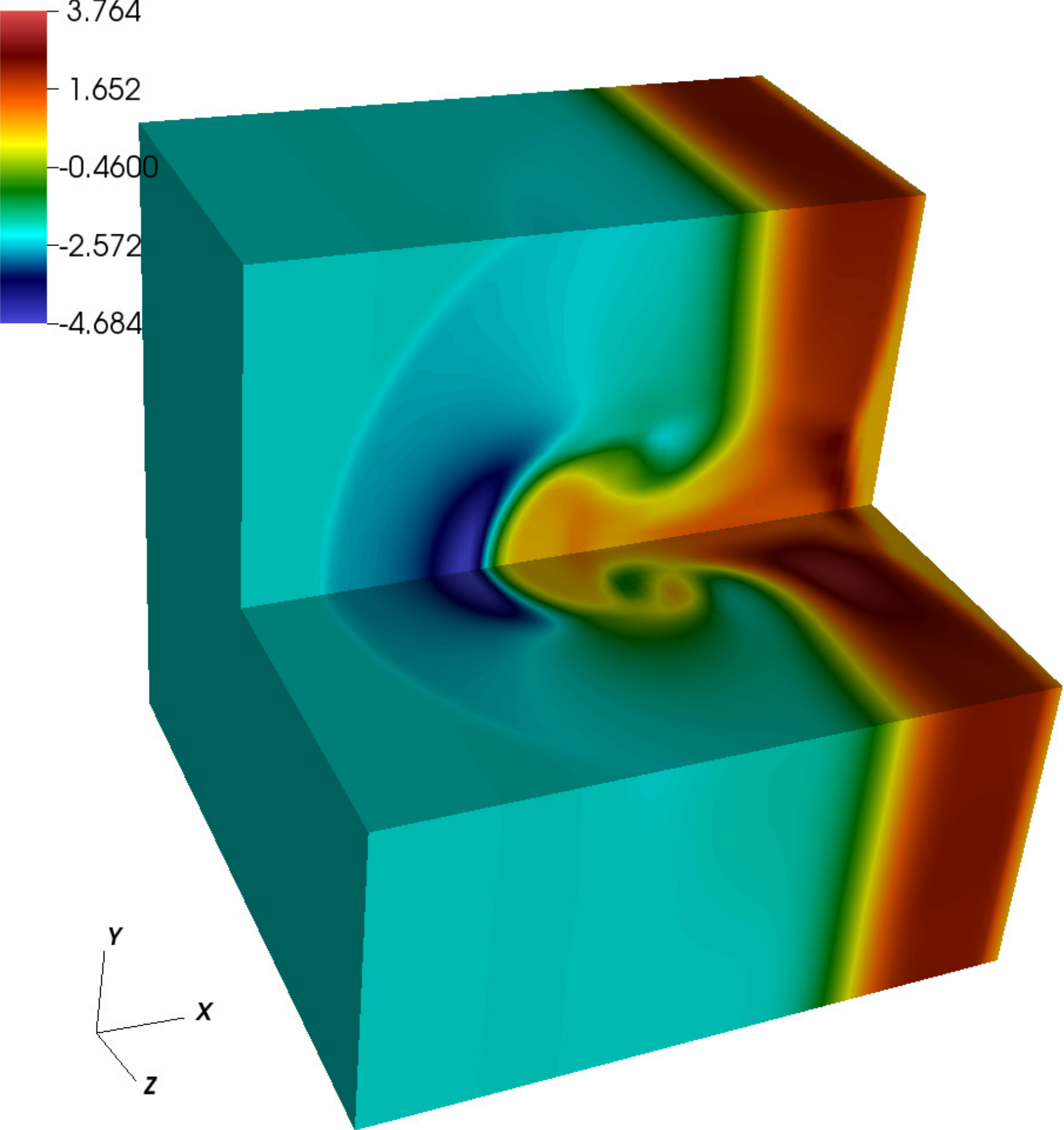}   \\
\end{tabular}
\end{center}
\label{fig:SC3DL61}
\end{figure}

To allow a quantitative comparison with the reference, we choose a coarser level, $L=6$, corresponding to $64^3$ cells. 
In Figure~\ref{fig:SC3DL61}, we present visualizations of the variables $\rho$ and $B_z$. Due to the dissipative effects and the refinement level, the structures of the solution are much smoother, compared to the ideal case. However, we can notice that the topology of the solution is maintained.
This adaptive simulation requires $76\%$ of the cells over time for $L=6$. When we refine the mesh, we obtain a reduction of the cells with an improvement of the CPU time. In particular, for a simulation with $L=7$ , $55\%$ of the cells are needed over time and there is a $35\%$ reduction in CPU time, which is significant in the computational context and reinforces the efficiency of the adaptive MR approach. The {\color{black} corresponding} errors are presented in Table~\ref{tab:sc3d-res}.

\begin{table}[H]
        \centering
        \caption{Errors obtained for the constant resistive 3D shock-cloud problem for the uniform and adaptive ($\epsilon^0=0.01$) cases, compared to the reference solution at level $L=6$.}
        \vspace{2mm}
        \small{
        \begin{tabular}{@{}clrr@{}}
        \toprule
        CARMEN--MHD &\multirow{2}{*}{Variables}       & \multicolumn{2}{c}{Errors $(\times 10^{-3})$} \\\cmidrule{3-4}
          solver &&   $\mathcal{L}^1$     &    $\mathcal{L}^2$    \\ 
           \cmidrule{1-4}
            \multirow{2}{*}{Uniform}
			& $\rho$ & 0.401 &  5.549 \\ 
   			& $B_z $ & 0.214 &  3.234 \\ 
        \toprule
            \multirow{2}{*}{Adaptive}
           	& $\rho$ & 0.403 &  5.558 \\ 
   			& $B_z $ & 0.214 &  3.230 \\ 
        \bottomrule
        \end{tabular}
        }      
        \label{tab:sc3d-res}
\end{table}

\subsection{Magnetic reconnection}
When we include Ohmic resistivity effects in the MHD equations, there is no conservation of magnetic flux anymore. 
This type of physical situation can change the topology of the magnetic field lines, allowing us to study different plasma problems, in particular, magnetic reconnection phenomena. 
The magnetic reconnection is a fundamental process in highly electrically conductive plasmas, which allows the conversion of magnetic energy to kinetic energy. 
It occurs when the magnetic field lines disconnect and reconnect again, changing its directions and restructuring the macroscopic plasma quantities.

We consider the Petschek reconnection model~\cite{petschek1964magnetic}, in which the reconnection rate is faster compared to the Sweet--Parker setting. The simulation proposed here was first presented in \cite{Jiang20121617} with the following configuration:
The initial condition for the magnetic reconnection is given by $\rho = 1$, $p = 0.1$, $\textbf{u}=\textbf{0}$, $B_x=0$ and
\begin{equation}
B_y \, = \, \begin{cases} \; -1,\text{ if } \; x<-0.05\\ \; \sin(\pi x/0.01),\text{ if } \; |x|\leq 0.05\\ \; 1,\text{ if} \; x>0.05\end{cases}
\end{equation}
The computational domain is $\Omega=[-0.5,0.5]\times[-2,2]$,
where the diffusion region is defined as $[-0.05,0.05]\times[0.2,0.2]$, and the resistivity inside this region is given  by \linebreak
$\eta(x,y) \, = \,  0.25\,\eta_0 \left( \cos(\pi x/0.1) + 1\right) \left(\cos(\pi y/0.4) + 1\right)$, 
where $\eta_0=0.00075$ is the initial resistivity. 

In absence of a reference solution, we compare our results to the ones presented in \cite{Jiang20121617}. The benchmark results are obtained on a $2048\times 4096$ high resolution mesh, while we use a $512\times 512$ mesh (corresponding to a refinement level $L=9$), employing a WENO scheme with Lax-Friedrich flux, a second order TVD Runge-Kutta time scheme and a damping approach at the boundaries.

In the presented simulations we use the final time $t=2.5$, the adiabatic constant $\gamma = 5/3$, the parabolic-hyperbolic correction parameter $\alpha_p = 0.4$, the Courant number $\nu = 0.4$, the threshold parameter $\epsilon = 0.0005$, and Neumann boundary conditions {\color{black} in all directions}.

\begin{figure}[H]
\psfrag{VY}{$u_y$}
\psfrag{MAG}{$B^2/2$}
\psfrag{x}{$x$}
\psfrag{GP}{p}
\psfrag{VX}{$\rho$}
\centering
\caption{Cuts in variables $p$, $u_y$, 
at $t=2.5$ and $L=9$, obtained with  non-constant resistive CARMEN--MHD with an  adaptive simulation using  $\epsilon_0=0.005$ (dotted);  and the reference solution presented in \cite[p. 1630, Fig. 11]{Jiang20121617}(solid line).} 
\begin{center}
\begin{tabular}{ccc}
$p$ &
$u_y$\\
\includegraphics[width=0.4\textwidth]{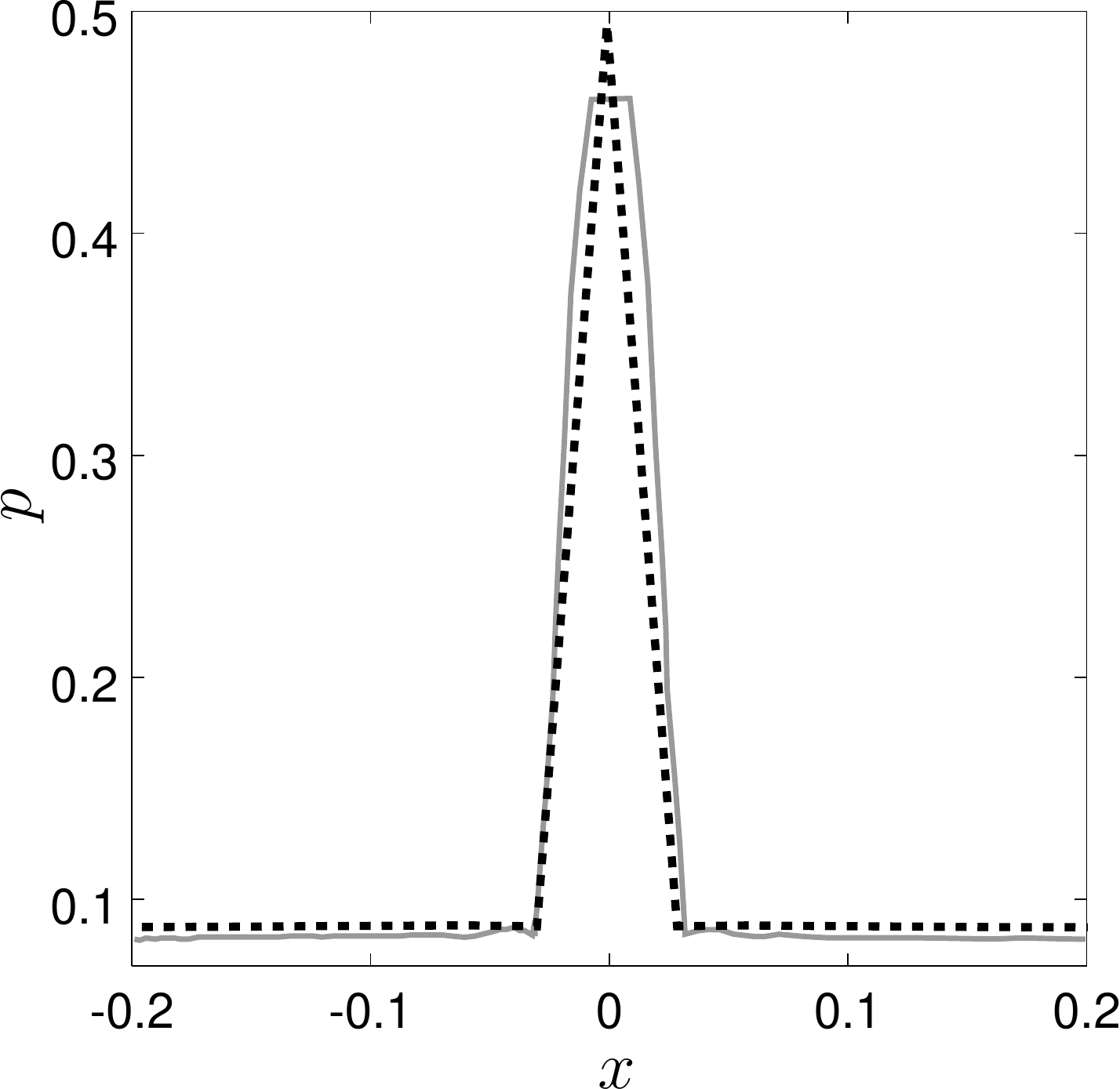} &
\includegraphics[width=0.4\textwidth]{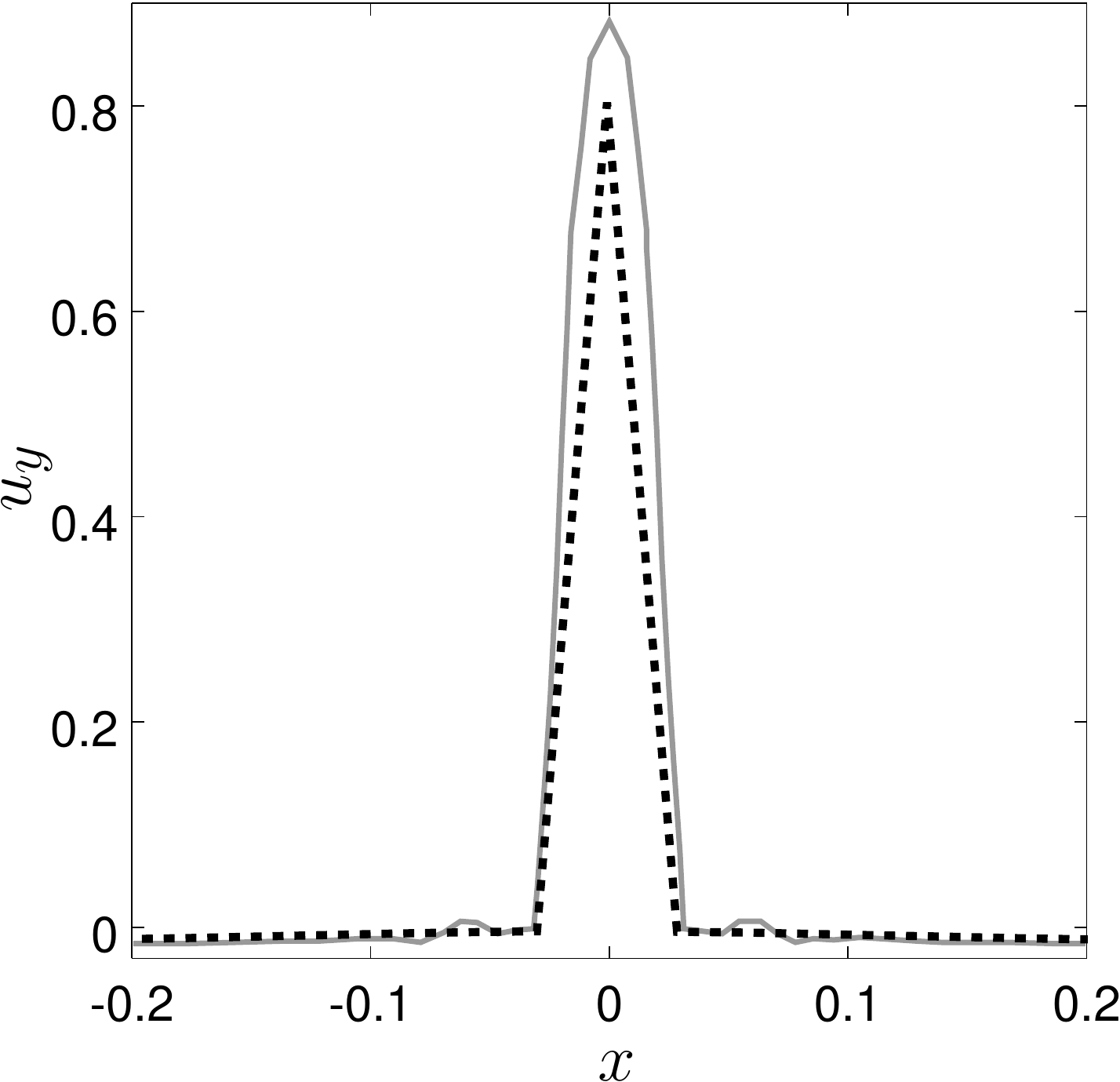} 
\end{tabular}
\end{center}
\label{fig:MR2DL9cut}
\end{figure}

The cuts at $y=0.5$ are presented for the variables $\rho$, $u_y$, 
in Figure~\ref{fig:MR2DL9cut}, for the interval $[-0.2,0,2]$. These cuts are similar between each other, however we can find some differences in the solutions, which can be attributed to the resolution or the chosen numerical scheme. High resolution can increase the accuracy of the numerical solution, however, the solution obtained at $L=9$ with the CARMEN--MHD code already presents the expected structures.

\begin{figure}[H]
\centering
\caption{Resistive magnetic reconnection problem: adaptive simulations with $\epsilon=0.005$ for the variables $B_x$,  $B_z$, and the adaptive mesh at $t=2.5$, and $L=9$ .}
\vspace{2mm}
\begin{center}
\begin{tabular}{ccc}
$B_x$ & $B_z$ & \hspace{7mm}adaptive mesh\\
\includegraphics[width=0.3\textwidth]{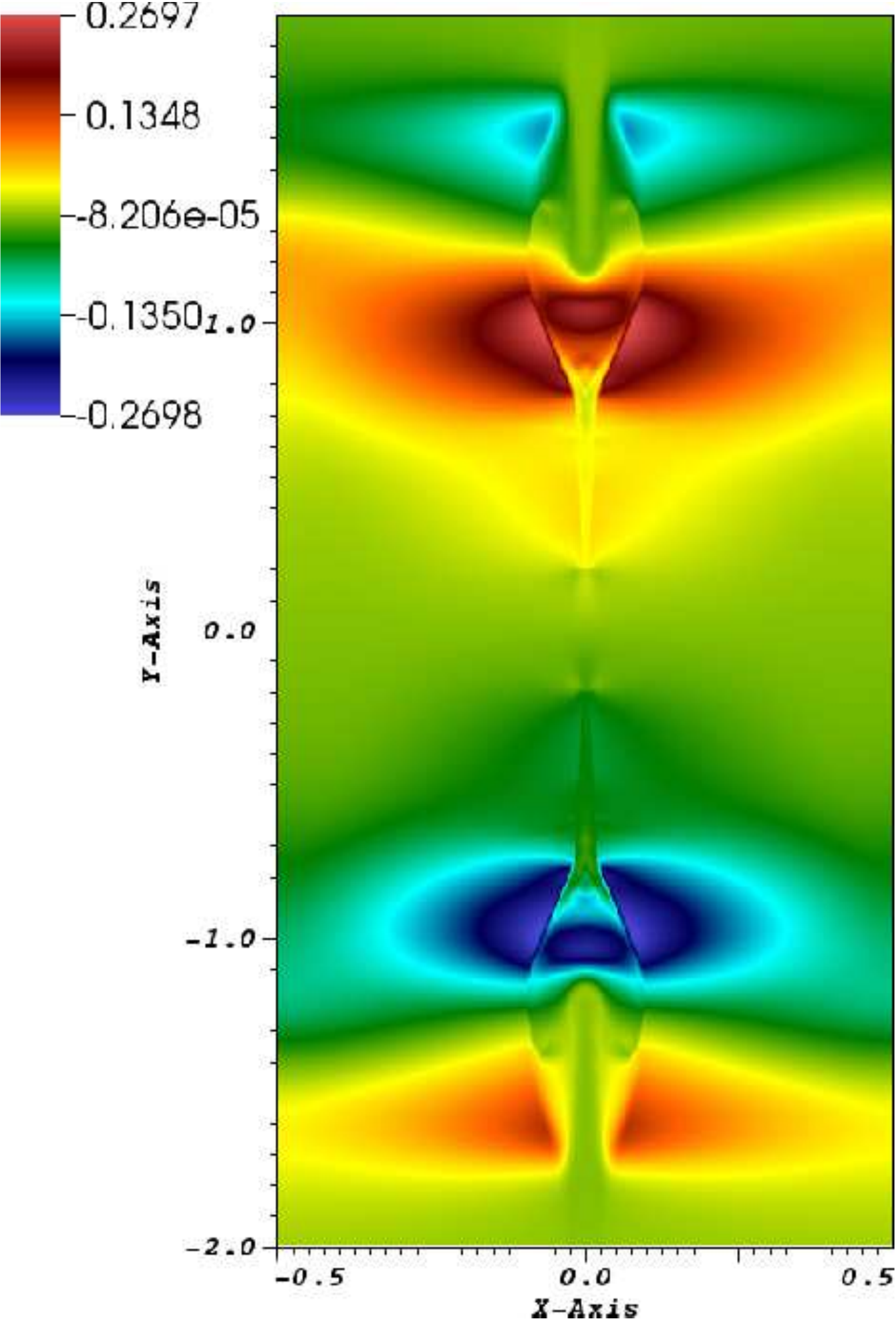} & \includegraphics[width=0.3\textwidth]{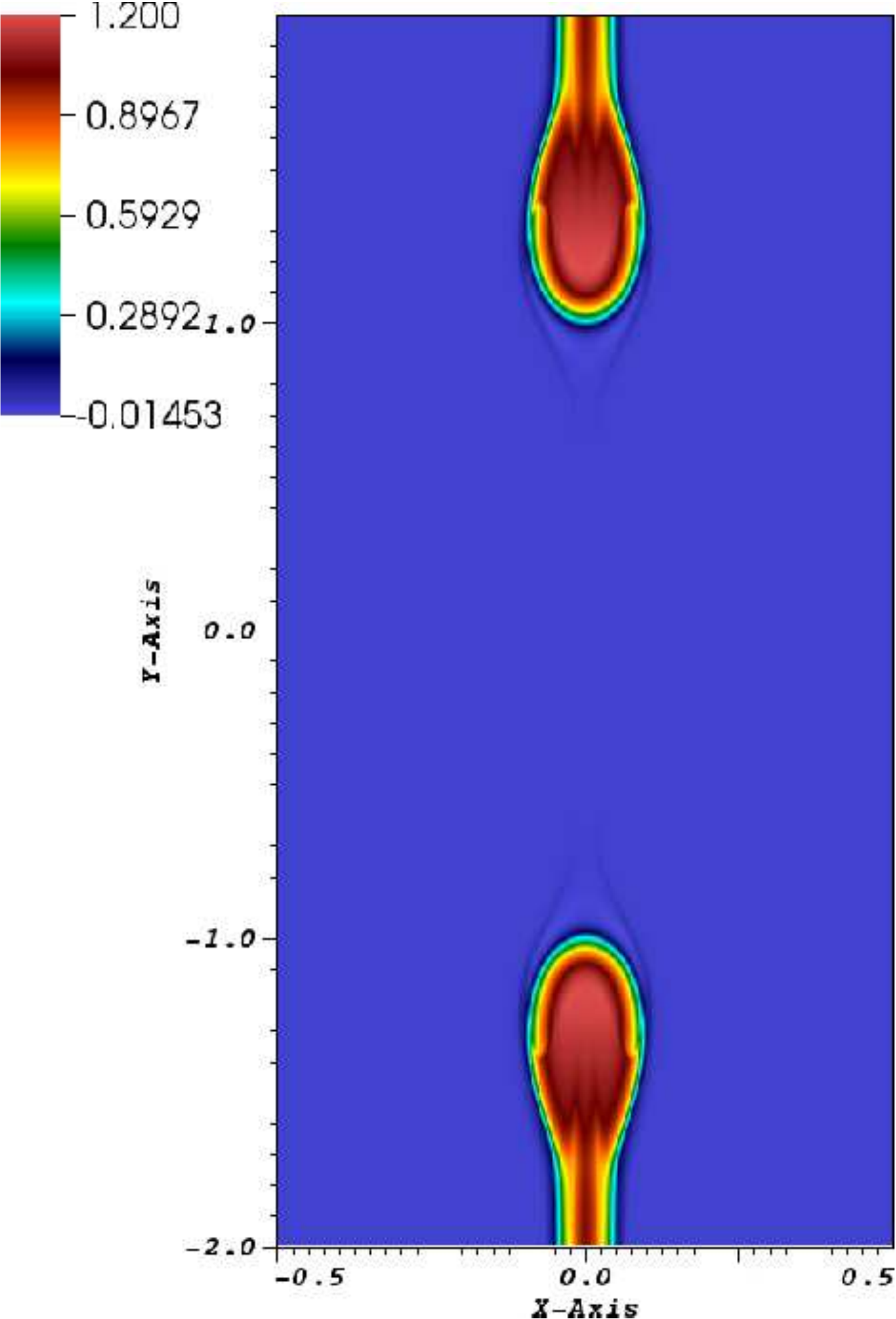} &
\includegraphics[width=0.3\textwidth]{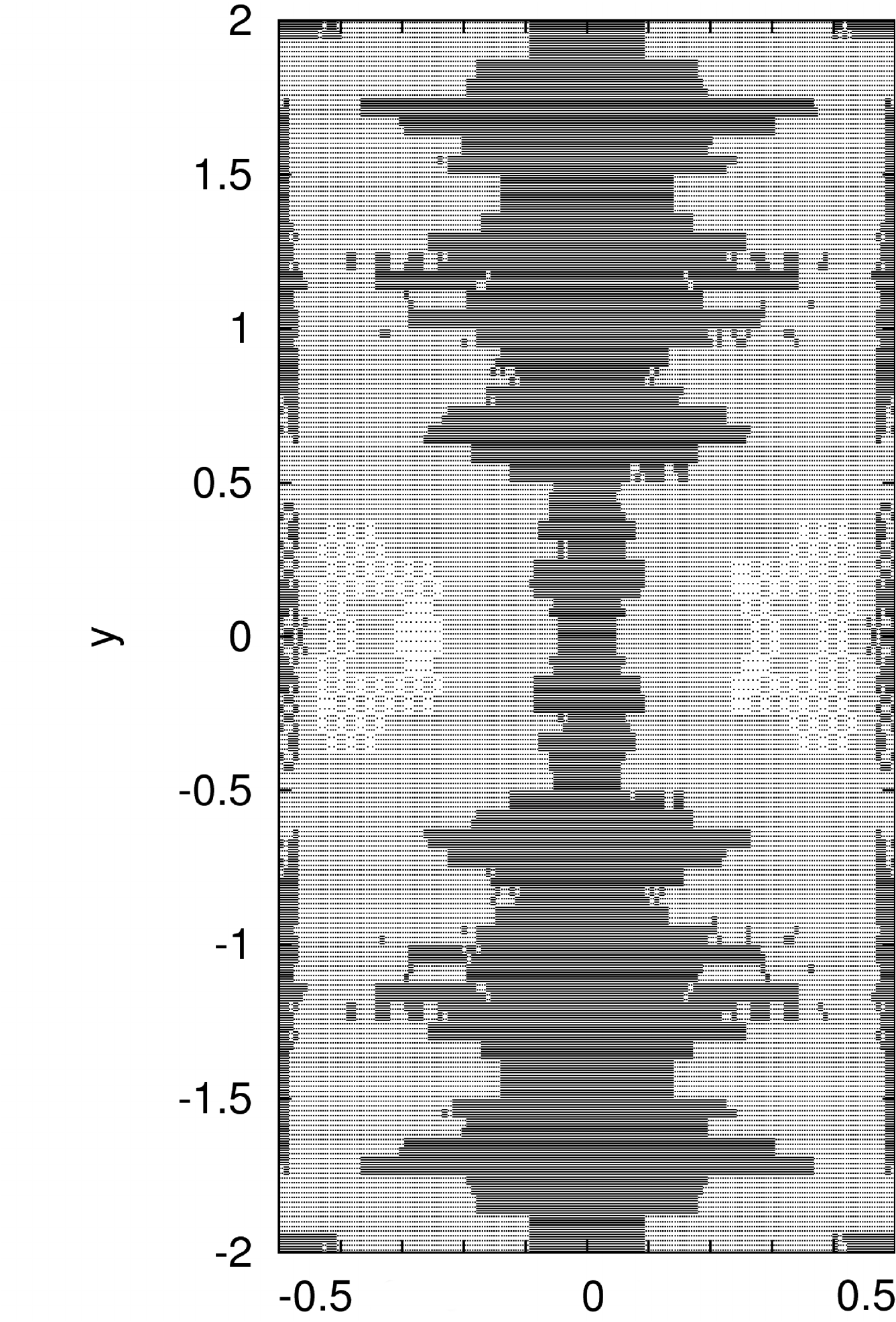} 
\end{tabular}
\end{center}
\label{fig:MR2DL9}
\end{figure}

The variables $B_x$ and $B_z$ obtained with the CARMEN--MHD code are shown in Figure~\ref{fig:MR2DL9}. 
The reconnection structure occurs well defined and is present in the computational domain in each variable. 
Neither numerical instabilities nor oscillations are observed during the simulation. 

The adaptive approach needs $55\%$ of the cells over time for this simulation, leading to a $40\%$ reduction in CPU time. We found that the threshold parameter $\epsilon=0.005$ is optimal for this case, since we can obtain a good {\color{black} compromise} between compression and coherent representation of the physical structures. When we choose a slightly larger value, e.g., $\epsilon=0.008$, the central structures of the problem are not well represented anymore. It is also possible to use smaller values for $\epsilon$, nevertheless this would lead to lower compression. 
{The majority of the refined cells in the adaptive mesh is located in the central region of the domain (corresponding to darker symbols). The other regions present a coarser refinement (corresponding to lighter symbols).}
 
 If we compare the {\color{black} structures present in the} variables $B_x$ and $B_y$ with the adaptive mesh, we can conclude that the mesh is efficiently adapted where {large discontinuities and diffusion regions are located}. This shows that the adaptive algorithm is indeed efficient to represent and identify automatically the structures of the problem. 
Moreover, we expect that higher resolution simulations {\color{black} will further improve the adaptive representation and thus the gains in memory and CPU time reduction}.

In this adaptive simulation the physical behavior of the magnetic reconnection phenomena is sustained.
In particular, we verified that the velocity field is compatible with the magnetic reconnection settings, in which the velocity direction initially points to the diffusion region and, after reconnection, its orientation does change.

\section{Conclusion}  
\label{sec4}

Fully adaptive numerical simulations using the CARMEN--MHD code were performed in {\color{black} two and three space dimensions} in order to verify the implementation and its computational efficiency. 
The numerical method is based on a finite volume discretization {\color{black} on Cartesian grids} and uses an adaptive multiresolution approach for introducing dynamically refined dyadic meshes.
{Our choice using Cartesian geometries instead of general geometries is motivated by the fact that multiresolution analysis is particularly attractive in this context. However generalization are possible in future work considering e.g. mutliresolution on triangles proposed by Cohen et al. \cite{Cohen2000} or more general tesselations borrowing techniques introduced in the field of computer graphics, e.g. using hierarchy refinement procedures~\cite{cohen2012}. This would allow designing adaptive multiresolution solvers on unstructured grids considering also complex geometries.}
Selected benchmarks were chosen in order to comprehend different physical and numerical {\color{black} phenomena}, and to ensure the correct behavior of the code in capturing the intrinsic topology of each situation.
The obtained results were then compared with the FLASH code which served as reference.

Quantitative and qualitative comparisons of the numerical solutions were {\color{black} carried out} and their convergence towards reference solutions was shown.
The obtained results are coherent for both, adaptive and uniform {\color{black} grid} approaches. The physical restrictions of the MHD model are maintained in the context of the numerical solution, contributing to the reliability of the results and its adequate reproduction.

Depending on the benchmark the structures present in the solution of the MHD model are located in different regions of the domain, varying according to the chosen variable. This type of situation is challenging in the context of the adaptivity, since the mesh must be adapted adequately. The adaptivity criteria used here were shown to be efficient for identifying the structures of the solution, even in cases where the structures do not present local features at all. We also showed that the numerical simulations are stable and do not require additional stabilization, e.g. adding numerical diffusion.

Moreover, we observed that it is possible to {\color{black} design} an optimal mesh adaptivity, by evaluating the relation between the threshold parameter and the approximation errors. The optimal choice can decrease the CPU time, {\color{black} while ensuring} the accuracy of the numerical solution. The adaptive multiresolution approach can increase significantly the computational gains of the simulations, even in non-parallel simulations. Thereby, this approach is shown to be computationally efficient to deal with the proposed MHD models.

Finally, let us mention that in the context of MHD, adaptive multiresolution computations, especially in 3D, are recent and in this work we presented {\color{black} their} potential, by showing their efficiency using adaptive meshes, while preserving {\color{black} the accuracy of the underlying discretization}.
We thus conclude that the verification of the CARMEN--MHD code was successful and we provide its open source-code and documentation for the community, in order to continue the research on MHD simulations of other interesting and challenging physical problems.

\section*{Acknowledgements}

The authors thank the FAPESP (Grant: $2015/ 25624-2$), 
CNPq (Grants: 
$302226/$ $2018-4, 307083/2017-9,
306038/2015-3, 302226/2018-4$), 
and FINEP (Grant: $0112052700$) for financial support of this research.
K.S. acknowledges partial support by the French Federation for Magnetic Fusion Studies (FR-FCM) and the Eurofusion consortium, funded by the Euratom research and training programme 2014-2018 and 2019-2020 under grant agreement No 633053. The views and opinions expressed herein do not necessarily reflect those of the European Commission. 
We are indebted to Eng. V. E. Menconi for his invaluable computational assistance.

\bibliographystyle{spmpsci}      
\bibliography{mhd}   

\end{document}